# Katyusha: The First Direct Acceleration of Stochastic Gradient Methods

(version 6)*

Zeyuan Allen-Zhu
zeyuan@csail.mit.edu
Princeton University / Institute for Advanced Study

March 18, 2016†

**Abstract**

Nesterov's momentum trick is famously known for accelerating gradient descent, and has been proven useful in building fast iterative algorithms. However, in the stochastic setting, counterexamples exist and prevent Nesterov's momentum from providing similar acceleration, even if the underlying problem is convex and finite-sum.

We introduce `Katyusha`, a direct, primal-only stochastic gradient method to fix this issue. In convex finite-sum stochastic optimization, `Katyusha` has an optimal accelerated convergence rate, and enjoys an optimal parallel linear speedup in the mini-batch setting.

The main ingredient is *Katyusha momentum*, a novel "negative momentum" on top of Nesterov's momentum. It can be incorporated into a variance-reduction based algorithm and speed it up, both in terms of *sequential and parallel* performance. Since variance reduction has been successfully applied to a growing list of practical problems, our paper suggests that in each of such cases, one could potentially try to give Katyusha a hug.

*We would like to specially thank Shai Shalev-Shwartz for useful feedbacks and suggestions on this paper, thank Blake Woodworth and Nati Srebro for pointer to their paper [56], thank Guanghui Lan for correcting our citation of [21], thank Weston Jackson, Xu Chen and Zhe Li for verifying the proofs and correcting typos, thank Hongzhou Lin for discussing the experiments, thank Robert Hannah and Lin Xiao for discussing lower bounds and parallel speedup, and thank anonymous reviewers for a number of writing suggestions. This paper is partially supported by an NSF Grant, no. CCF-1412958, and a Microsoft Research Grant, no. 0518584. Any opinions, findings and conclusions or recommendations expressed in this material are those of the author(s) and do not necessarily reflect the views of NSF or Microsoft.

†V1 of this paper appeared on this date. V2 includes experiments, V3/V4 polish writing, and V5/V6 is the journal revision and includes missing proofs for some extended settings. This paper has appeared in Symposium on Theory of Computing (STOC 2017) and Journal of Machine Learning Research (JMLR 2018).

# 1 Introduction

In large-scale machine learning, the number of data examples is usually very large. To search for the optimal solution, one often uses *stochastic gradient methods* which only require one (or a small batch of) random example(s) per iteration in order to form an *estimator* of the full gradient.

While full-gradient based methods can enjoy an *accelerated* (and optimal) convergence rate if Nesterov's momentum trick is used [42–44], theory for stochastic gradient methods are generally lagging behind and less is known for their acceleration.

At a high level, momentum is *dangerous* if stochastic gradients are present. If some gradient estimator is very inaccurate, then adding it to the momentum and moving further in this direction (for every future iteration) may hurt the convergence performance. In other words, when naively equipped with momentum, stochastic gradient methods are "very prune to error accumulation" [31] and do *not* yield accelerated convergence rates in general.[1]

In this paper, we show that at least for convex optimization purposes, such an issue can be solved with a novel "negative momentum" that can be added on top of momentum. We obtain accelerated and the first optimal convergence rates for stochastic gradient methods. As one of our side results, under this "negative momenmtum," our new method enjoys a linear speedup in the parallel (i.e., mini-match) setting. We hope our new insight could potentially deepen our understanding to the theory of accelerated methods.

**Problem Definition.** Consider the following composite convex minimization problem

$$\min_{x\in\mathbb{R}^d}\Big\{F(x)\stackrel{\text{def}}{=}f(x)+\psi(x)\stackrel{\text{def}}{=}\frac{1}{n}\sum_{i=1}^{n}f_i(x)+\psi(x)\Big\}\ . \tag{1.1}$$

Here, $f(x)=\frac{1}{n}\sum_{i=1}^{n}f_i(x)$ is a convex function that is a finite average of $n$ convex, smooth functions $f_i(x)$, and $\psi(x)$ is convex, lower semicontinuous (but possibly non-differentiable) function, sometimes referred to as the *proximal* function. We mostly focus on the case when $\psi(x)$ is $\sigma$-strongly convex and each $f_i(x)$ is $L$-smooth. (Both these assumptions can be removed and we shall discuss that later.) We look for approximate minimizers $x\in\mathbb{R}^d$ satisfying $F(x)\le F(x^*)+\varepsilon$, where $x^*\in\arg\min_x\{F(x)\}$.

Problem (1.1) arises in many places in machine learning, statistics, and operations research. All convex *regularized empirical risk minimization (ERM)* problems such as Lasso, SVM, Logistic Regression, fall into this category (see Section 1.3). Efficient stochastic methods for Problem (1.1) have also inspired stochastic algorithms for neural nets [5, 29, 34] as well as SVD, PCA, and CCA [7, 9, 26].

We summarize the history of stochastic gradient methods for Problem (1.1) in three eras.

**The First Era: Stochastic Gradient Descent (SGD).**

Recall that stochastic gradient methods iteratively perform the following update

$$\text{stochastic gradient iteration:}\qquad x_{k+1}\leftarrow\arg\min_{y\in\mathbb{R}^d}\Big\{\frac{1}{2\eta}\|y-x_k\|_2^2+\langle\widetilde{\nabla}_k,y\rangle+\psi(y)\Big\}\ ,$$

where $\eta$ is the step length and $\widetilde{\nabla}_k$ is a random vector satisfying $\mathbb{E}[\widetilde{\nabla}_k]=\nabla f(x_k)$ and is referred to as the *gradient estimator*. If the proximal function $\psi(y)$ equals zero, the update reduces to $x_{k+1}\leftarrow x_k-\eta\widetilde{\nabla}_k$. A popular choice for the gradient estimator is to set $\widetilde{\nabla}_k=\nabla f_i(x_k)$ for some random index $i \in [n]$ per iteration, and methods based on this choice are known as *stochastic gradient descent (SGD)* [18, 59]. Since computing $\nabla f_i(x)$ is usually $n$ times faster than that of $\nabla f(x)$, SGD enjoys a low per-iteration cost as compared to full-gradient methods; however, SGD cannot converge at a rate faster than $1/\varepsilon$ even if $F(\cdot)$ is strongly convex and smooth.

[1] In practice, experimentalists have observed that momentums could sometimes help if stochastic gradient iterations are used. However, the so-obtained methods (1) sometimes fail to converge in an accelerated rate, (2) become unstable and hard to tune, and (3) have no support theory behind them. See Section 7.1 for an experiment illustrating that, even for convex stochastic optimization.

**The Second Era: Variance Reduction Gives Faster Convergence.**

The convergence rate of SGD can be further improved with the so-called *variance-reduction* technique, first proposed by Schmidt et al. [49] (solving a sub-case of Problem (1.1)) and then followed by many others [13, 22, 29, 39, 40, 51–53, 57, 58]. In these cited results, the authors have shown that SGD converges much faster if one makes a better choice of the gradient estimator $\widetilde{\nabla}_k$ so that its variance reduces as $k$ increases. One way to choose this estimator can be described as follows [29, 58]. Keep a *snapshot* vector $\widetilde{x} = x_k$ that is updated once every $m$ iterations (where $m$ is some parameter usually around $2n$), and compute the full gradient $\nabla f(\widetilde{x})$ only for such snapshots. Then, set

$$\widetilde{\nabla}_k = \nabla f_i(x_k) - \nabla f_i(\widetilde{x}) + \nabla f(\widetilde{x}) \ . \tag{1.2}$$

This choice of gradient estimator ensures that its variance approaches to zero as $k$ grows. Furthermore, the number of stochastic gradients (i.e., the number of computations of $\nabla f_i(x)$ for some $i$) required to reach an $\varepsilon$-approximate minimizer of Problem (1.1) is only $O\big((n + \frac{L}{\sigma}) \log \frac{1}{\varepsilon}\big)$. Since it is often denoted by $\kappa \stackrel{\text{def}}{=} L/\sigma$ the condition number of the problem, we rewrite the above iteration complexity as $O\big((n+\kappa)\log\frac{1}{\varepsilon}\big)$.

Unfortunately, the iteration complexities of all known variance-reduction based methods have a linear dependence on $\kappa$. It was an open question regarding how to obtain an *accelerated* stochastic gradient method with an optimal $\sqrt{\kappa}$ dependency.

**The Third Era: Acceleration Gives Fastest Convergence.**

This open question was partially solved recently by the APPA [24] and Catalyst [36] reductions, both based on an outer-inner loop structure first proposed by Shalev-Shwartz and Zhang [54]. We refer to both of them as Catalyst in this paper. Catalyst solves Problem (1.1) using $O\big((n + \sqrt{n\kappa}) \log \kappa \log \frac{1}{\varepsilon}\big)$ stochastic gradient iterations, through a logarithmic number of calls to a variance-reduction method.[2] However, Catalyst is still imperfect for the following reasons:

- Optimality. Catalyst does not match the optimal $\sqrt{\kappa}$ dependence [56] and has an extra $\log \kappa$ factor. It yields suboptimal rate $\frac{\log^4 T}{T^2}$ if the objective is not strongly convex or is non-smooth; and it yields suboptimal rate $\frac{\log^4 T}{T}$ if the objective is both non-strongly convex and non-smooth.[3]
- Practicality. To the best of our knowledge, Catalyst is not very practical since each of its inner iterations needs to be very accurately executed. This makes the stopping criterion hard to be tuned, and makes Catalyst sometimes run slower than non-accelerated variance-reduction methods. We have also confirmed this in our experiments.
- Parallelism. To the best of our knowledge, Catalyst does not give competent parallel performance (see Section 1.2). If $b \in \{1, \dots, n\}$ stochastic gradients (instead of one) are computed in each iteration, the number of iterations of Catalyst reduces by $O(\sqrt{b})$. In contrast, the best parallel speedup one can hope for is "linear speedup": that is, to reduce the number of iterations by a factor of $O(b)$ for $b \le \sqrt{n}$.
- Generality. To the best of our knowledge, being a reduction-based method, Catalyst does not support non-Euclidean norm smoothness (see Section 1.2).

[2]Note that $n + \sqrt{n\kappa}$ is always less than $O(n+\kappa)$.

[3]Obtaining *optimal* rates is one of the main goals in optimization and machine learning. For instance, obtaining the optimal $1/T$ rate for online learning was a very meaningful result, even though the $\log T/T$ rate was known [27, 48].

Another acceleration method by Lan and Zhou [33] is based on a primal-dual analysis that also has suboptimal convergence rates and suboptimal parallel speedup like Catalyst. Their method requires $n$ times more storage compared with Catalyst for solving Problem (1.1).

In sum, it is desirable and also an open question to develop a *direct*, *primal-only*, and *optimal* accelerated stochastic gradient method without using reductions. This could have both theoretical and practical impacts to the problems that fall into the general framework of (1.1), and potentially deepen our understanding to acceleration in stochastic settings.

## 1.1 Our Main Results and High-Level Ideas

We develop a direct, accelerated stochastic gradient method `Katyusha` for Problem (1.1) in

$$\boxed{O\big((n+\sqrt{n\kappa})\log(1/\varepsilon)\big) \text{ stochastic gradient iterations (see Theorem 2.1).}}$$

This gives both optimal dependency on $\kappa$ and on $\varepsilon$ which was not obtained before for stochastic gradient methods. In addition, if $F(\cdot)$ is non-strongly convex (non-SC), `Katyusha` converges to an $\varepsilon$-minimizer in

$$\boxed{O\big(n\log(1/\varepsilon)+\sqrt{nL/\varepsilon}\,\big) \text{ stochastic gradient iterations (see Corollary 3.7).}}$$

This gives an optimal $\varepsilon \propto \frac{n}{T^2}$ rate where in contrast Catalyst has rate $\varepsilon \propto \frac{n\log^4 T}{T^2}$. The lower bound from Woodworth and Srebro [56] is $\Omega\big(n+\sqrt{nL/\varepsilon}\big)$.

**Our Algorithm.** If ignoring the proximal term $\psi(\cdot)$ and viewing it as zero, our `Katyusha` method iteratively perform the following updates for $k = 0, 1, \dots$ :

- $x_{k+1} \leftarrow \tau_1 z_k + \tau_2 \widetilde{x} + (1-\tau_1-\tau_2) y_k$; (so $x_{k+1} = y_k + \tau_1(z_k - y_k) + \tau_2(\widetilde{x} - y_k)$)
- $\widetilde{\nabla}_{k+1} \leftarrow \nabla f(\widetilde{x}) + \nabla f_i(x_{k+1}) - \nabla f_i(\widetilde{x})$ where $i$ is a random index in $[n]$;
- $y_{k+1} \leftarrow x_{k+1} - \frac{1}{3L}\widetilde{\nabla}_{k+1}$, and
- $z_{k+1} \leftarrow z_k - \alpha\widetilde{\nabla}_{k+1}$.

Above, $\widetilde{x}$ is a snapshot point which is updated every $m$ iterations, $\widetilde{\nabla}_{k+1}$ is the gradient estimator defined in the same way as (1.2), $\tau_1, \tau_2 \in [0,1]$ are two momentum parameters, and $\alpha$ is a parameter that is equal to $\frac{1}{3\tau_1 L}$. The reason for keeping three vector sequences $(x_k, y_k, z_k)$ is a common ingredient that can be found in all existing accelerated methods.[4]

**Our New Technique– Katyusha Momentum.** The most interesting ingredient of `Katyusha` is the novel choice of $x_{k+1}$ which is a convex combination of $y_k$, $z_k$, and $\widetilde{x}$. Our theory suggests the parameter choices $\tau_2 = 0.5$ and $\tau_1 = \min\{\sqrt{n\sigma/L}, 0.5\}$ and they work well in practice too. To explain this novel combination, let us recall the classical "momentum" view of accelerated methods.

In a classical accelerated gradient method, $x_{k+1}$ is only a convex combination of $y_k$ and $z_k$ (or equivalently, $\tau_2 = 0$ in our formulation). At a high level, $z_k$ plays the role of "momentum" which adds a weighted sum of the gradient history into $y_{k+1}$. As an illustrative example, suppose $\tau_2 = 0$, $\tau_1 = \tau$, and $x_0 = y_0 = z_0$. Then, one can compute that

$$y_k = \begin{cases} x_0 - \frac{1}{3L}\widetilde{\nabla}_1, & k=1; \\ x_0 - \frac{1}{3L}\widetilde{\nabla}_2 - \big((1-\tau)\frac{1}{3L} + \tau\alpha\big)\widetilde{\nabla}_1, & k=2; \\ x_0 - \frac{1}{3L}\widetilde{\nabla}_3 - \big((1-\tau)\frac{1}{3L} + \tau\alpha\big)\widetilde{\nabla}_2 - \big((1-\tau)^2\frac{1}{3L} + (1-(1-\tau)^2)\alpha\big)\widetilde{\nabla}_1, & k=3. \end{cases}$$

[4] One can of course rewrite the algorithm and keep track of only two vectors per iteration during implementation. This will make the algorithm statement less clean so we refrain from doing so in this paper.

Since $\alpha$ is usually much larger than $1/3L$, the above recursion suggests that the contribution of a fixed gradient $\widetilde{\nabla}_t$ gradually increases as time goes. For instance, the weight on $\widetilde{\nabla}_1$ is increasing because $\frac{1}{3L} < \big((1-\tau)\frac{1}{3L} + \tau\alpha\big) < \big((1-\tau)^2\frac{1}{3L} + (1-(1-\tau)^2)\alpha\big)$ . This is known as "momentum" which is at the heart of all accelerated first-order methods.

In `Katyusha`, we put a "magnet" around $\widetilde{x}$, where we choose $\widetilde{x}$ to be essentially "the average $x_t$ of the most recent $n$ iterations". Whenever we compute the next $x_{k+1}$, it will be attracted by the magnet $\widetilde{x}$ with weight $\tau_2 = 0.5$. This is a strong magnet: it ensures that $x_{k+1}$ is not too far away from $\widetilde{x}$ so the gradient estimator remains "accurate enough". This can be viewed as a "negative momentum" component, because the magnet retracts $x_{k+1}$ back to $\widetilde{x}$ and this can be understood as "counteracting a fraction of the positive momentum incurred from earlier iterations."

*We call it the Katyusha momentum.*

This summarizes the high-level idea behind our `Katyusha` method. We remark here if $\tau_1 = \tau_2 = 0$, `Katyusha` becomes almost identical to SVRG [29, 58] which is a variance-reduction based method.

## 1.2 Our Side Results

**Parallelism / Mini-batch.** Instead of using a single $\nabla f_i(\cdot)$ per iteration, for any stochastic gradient method, one can replace it with the average of $b$ stochastic gradients $\frac{1}{b}\sum_{i\in S}\nabla f_i(\cdot)$, where $S$ is a random subset of $[n]$ with cardinality $b$. This is known as the *mini-batch* technique and it allows the stochastic gradients to be computed in a distributed manner, using up to $b$ processors.

Our `Katyusha` method trivially extends to this mini-batch setting. For instance, at least for $b \in \{1, 2, \dots, \lceil\sqrt{n}\rceil\}$, `Katyusha` enjoys a *linear speedup* in the parallel running time. In other words, if ignoring communication overhead,

`Katyusha` can be distributed to $b \leq \sqrt{n}$ machines with a parallel speed-up factor $b$.

In contrast, to the best of our knowledge, without any additional assumption, (1) non-accelerated methods such as SVRG or SAGA are not known to enjoy any parallel speed-up; (2) Catalyst enjoys a parallel speed-up factor of only $\sqrt{b}$. Details are in Section 5.

**Non-Uniform Smoothness.** If each $f_i(\cdot)$ has a possibly different smooth parameter $L_i$ and $\overline{L} = \frac{1}{n}\sum_{i=1}^n L_i$, then an naive implementation of `Katyusha` only gives a complexity that depends on $\max_i L_i$ but not $\overline{L}$. In such a case, we can select the random index $i \in [n]$ with probability proportional to $L_i$ per iteration to slightly improve the total running time.

Furthermore, suppose $f(x) = \frac{1}{n}\sum_{i=1}^n f_i(x)$ is smooth with parameter $L$, it satisfies $\overline{L} \in [L, nL]$. One can ask whether or not $L$ influences the performance of `Katyusha`. We show that, in the mini-batch setting when $b$ is large, the total complexity becomes a function on $L$ as opposed to $\overline{L}$. The details are in Section 5.

**A Precise Statement.** Taking into account both the mini-batch parameter $b$ and the non-uniform smoothness parameters $L$ and $\overline{L}$, we show `Katyusha` solves Problem (1.1) in

$O\Big(\big(n + b\sqrt{L/\sigma} + \sqrt{n\overline{L}/\sigma}\big) \cdot \log\frac{1}{\varepsilon}\Big)$ stochastic gradient computations (see Theorem 5.4)

**Non-Euclidean Norms.** If the smoothness of each $f_i(x)$ is with respect to a non-Euclidean norm (such as the well known $\ell_1$ norm case over the simplex), our main result still holds. Our update on the $y_{k+1}$ side becomes the non-Euclidean norm gradient descent, and our update on the $z_{k+1}$ side becomes the non-Euclidean norm mirror descent. We include such details in Section 6. To the best of our knowledge, most known accelerated methods (including Catalyst, AccSDCA and APCG) do not work with non-Euclidean norms. SPDC can be revised to work with non-Euclidean norms, see [15].

**Remark on Katyusha Momentum Weight $\tau_2$.** To provide the simplest proof, we choose $\tau_2 = 1/2$ which also works well in practice. Our proof trivially generalizes to all constant values $\tau_2 \in (0,1)$, and it could be beneficial to tune $\tau_2$ for different datasets. However, for a stronger comparison, in our experiments we refrain from tuning $\tau_2$: by fixing $\tau_2 = 1/2$ and without increasing parameter tuning difficulties, `Katyusha` already outperforms most of the state-of-the-arts.

In the mini-batch setting, it turns out the best theoretical choice is essentially $\tau_2 = \frac{1}{2b}$, where $b$ is the size of the mini-batch. In other words, the larger the mini-batch size, the smaller weight we want to give to Katyusha momentum. This should be intuitive, because when $b = n$ we are almost in the full-gradient setting and do not need Katyusha momentum.

## 1.3 Applications: Optimal Rates for Empirical Risk Minimization

Suppose we are given $n$ feature vectors $a_1, \dots, a_n \in \mathbb{R}^d$ corresponding to $n$ data samples. Then, the *empirical risk minimization (ERM)* problem is to study Problem (1.1) when each $f_i(x)$ is "rank-one" structured: $f_i(x) \stackrel{\text{def}}{=} g_i(\langle a_i, x\rangle)$ for some loss function $g_i \colon \mathbb{R} \to \mathbb{R}$. Slightly abusing notation, we write $f_i(x) = f_i(\langle a_i, x\rangle)$.[5] In such a case, Problem (1.1) becomes

$$\text{ERM:} \quad \min_{x\in\mathbb{R}^d} \Big\{ F(x) \stackrel{\text{def}}{=} f(x) + \psi(x) \stackrel{\text{def}}{=} \frac{1}{n}\sum_{i=1}^n f_i(\langle a_i, x\rangle) + \psi(x) \Big\} \ . \tag{1.3}$$

Without loss of generality, we assume each $a_i$ has norm 1 because otherwise one can scale $f_i(\cdot)$ accordingly. As summarized for instance in [6], there are four interesting cases of ERM problems, all can be written in the form of (1.3):

Case 1: $\psi(x)$ is $\sigma$-SC and $f_i(x)$ is $L$-smooth. Examples: ridge regression, elastic net;

Case 2: $\psi(x)$ is non-SC and $f_i(x)$ is $L$-smooth. Examples: Lasso, logistic regression;

Case 3: $\psi(x)$ is $\sigma$-SC and $f_i(x)$ is non-smooth. Examples: support vector machine;

Case 4: $\psi(x)$ is non-SC and $f_i(x)$ is non-smooth. Examples: $\ell_1$-SVM.

**Known Results.** For all of the four ERM cases above, accelerated stochastic methods were introduced in the literature, most notably AccSDCA [54], APCG [37], SPDC [60]. All these methods rely on dual analysis so have suboptimal convergence rates for Cases 2, 3 and 4. (In fact, they also have the suboptimal dependence on the condition number $L/\sigma$ for Case 1.)

*Remark* 1.1. The best known rate was $\frac{\log(1/\varepsilon)}{\sqrt{\varepsilon}}$, $\frac{\log(1/\varepsilon)}{\sqrt{\varepsilon}}$, or $\frac{\log(1/\varepsilon)}{\varepsilon}$ respectively for Cases 2, 3, or 4, and is a factor $\log(1/\varepsilon)$ worse than optimal [56]. It is open to design a stochastic gradient method to match the optimal rates. In particular, Dang and Lan [21] provided an interesting attempt to remove such log factors but using a non-classical notion of convergence.[6]

Besides the log factor loss in the running time,[7] the aforementioned methods suffer from several other issues that most dual-based methods also suffer. First, they only apply to ERM problems but not to the more general Problem (1.1). Second, they require proximal updates with respect to the Fenchel conjugate $f_i^*(\cdot)$ which is sometimes unpleasant to work with. Third, their performances cannot benefit from the implicit strong convexity in $f(\cdot)$. All of these issues together make these methods sometimes even outperformed by primal-only non-accelerated ones, such as SAGA [22] or SVRG [29, 58].

[5] Assuming "rank-one" simplifies the notations; all of the results stated in this subsection generalize to constant-rank structured functions $f_i(x)$.

[6] Dang and Lan [21] work in a primal-dual $\phi(x,y)$ formulation of Problem (1.1), and produce a primal-dual pair $(x,y)$ so that for every fixed $(u,v)$, the expectation $\mathbb{E}[\phi(x,v) - \phi(u,y)] \le \varepsilon$. Unfortunately, to ensure $x$ is an $\varepsilon$-approximate minimizer of Problem (1.1), one needs the stronger $\mathbb{E}[\max_{(u,v)} \phi(x,v) - \phi(u,y)] \le \varepsilon$ to hold.

[7] In fact, dual-based methods have to suffer from a log factor loss in the convergence rate. This is so because even for Case 1 of Problem (1.3), converting an $\varepsilon$-maximizer for the dual objective to the primal, one only obtains an $n\kappa\varepsilon$-minimizer on the primal objective. As a result, algorithms like APCG who directly work on the dual, algorithms like SPDC who maintain both primal and dual variables, and algorithms like RPDG [33] that are primal-like but still use dual analysis, have to suffer from a log loss in the convergence rates.

**Our Results.** `Katyusha` simultaneously closes the gap for all of the three classes of problems with the help from the optimal reductions developed in [6]. We obtain an $\varepsilon$-approximate minimizer for Case 2 in $O\big(n\log\frac{1}{\varepsilon}+\frac{\sqrt{nL}}{\sqrt{\varepsilon}}\big)$ iterations, for Case 3 in $O\big(n\log\frac{1}{\varepsilon}+\frac{\sqrt{n}}{\sqrt{\sigma\varepsilon}}\big)$ iterations, and for Case 4 in $O\big(n\log\frac{1}{\varepsilon}+\frac{\sqrt{n}}{\varepsilon}\big)$ iterations. None of the existing accelerated methods can lead to such optimal rates even if the optimal reductions are used.

Woodworth and Srebro [56] proved the tightness of our results. They showed lower bounds $\Omega\big(n+\frac{\sqrt{nL}}{\sqrt{\varepsilon}}\big)$, $\Omega\big(n+\frac{\sqrt{n}}{\sqrt{\sigma\varepsilon}}\big)$, and $\Omega\big(n+\frac{\sqrt{n}}{\varepsilon}\big)$ for Cases 2, 3, and 4 respectively.[8]

### 1.4 Roadmap

- In Section 2, we state and prove our theorem on `Katyusha` for the strongly convex case.
- In Section 3, we apply `Katyusha` to non-strongly convex or non-smooth cases by reductions.
- In Section 4, we provide a *direct* algorithm $\texttt{Katyusha}^{\textsf{ns}}$ for the non-strongly case.
- In Section 5, we generalize `Katyusha` to mini-batch and non-uniform smoothness.
- In Section 6, we generalize `Katyusha` to the non-Euclidean norm setting.
- In Section 7, we provide an empirical evaluation to illustrate the necessity of Katyusha momentum, and the practical performance of `Katyusha`.

### 1.5 Notations

Throughout this paper (except Section 6), we denote by $\|\cdot\|$ the Euclidean norm. We denote by $\nabla f(x)$ the full gradient of function $f$ if it is differentiable, or any of its subgradients if $f$ is only Lipschitz continuous. Recall some classical definitions on strong convexity (SC) and smoothness.

**Definition 1.2.** *For a convex function $f\colon \mathbb{R}^n \to \mathbb{R}$,*
- *$f$ is $\sigma$-strongly convex if $\forall x, y \in \mathbb{R}^n$, it satisfies $f(y) \geq f(x) + \langle \nabla f(x), y-x\rangle + \frac{\sigma}{2}\|x-y\|^2$.*
- *$f$ is $L$-smooth if $\forall x, y \in \mathbb{R}^n$, it satisfies $\|\nabla f(x) - \nabla f(y)\| \leq L\|x-y\|$.*

## 2 Katyusha in the Strongly Convex Setting

We formally introduce our `Katyusha` algorithm in Algorithm 1. It follows from our high-level description in Section 1.1, and we make several remarks here behind our specific design.

- `Katyusha` *is divided into epochs each consisting of $m$ iterations. In theory, $m$ can be anything linear in $n$. We let snapshot $\widetilde{x}$ be a weighted average of $y_k$ in the most recent epoch.*

  $\widetilde{x}$ and $\widetilde{\nabla}_k$ correspond to a standard design on variance-reduced gradient estimators, called SVRG [29, 58]. The practical recommendation is $m = 2n$ [29]. Our choice $\widetilde{\nabla}_k$ is independent from our acceleration techniques, and we expect our result continues to apply to other

[8]More precisely, their lower bounds for Cases 3 and 4 are $\Omega\big(\min\big\{\frac{1}{\sigma\varepsilon}, n+\frac{\sqrt{n}}{\sqrt{\sigma\varepsilon}}\big\}\big)$ and $\Omega\big(\min\big\{\frac{1}{\varepsilon^2}, n+\frac{\sqrt{n}}{\varepsilon}\big\}\big)$. However, since the vanilla SGD requires $O(\frac{1}{\sigma\varepsilon})$ and $O(\frac{1}{\varepsilon^2})$ iterations for Cases 3 and 4, such lower bounds are matched by combining the best between `Katyusha` and SGD.

**Algorithm 1** `Katyusha`$(x_0, S, \sigma, L)$

$m \leftarrow 2n$; ◇ *epoch length*
$\tau_2 \leftarrow \frac{1}{2}$, $\tau_1 \leftarrow \min\left\{\frac{\sqrt{m\sigma}}{\sqrt{3L}}, \frac{1}{2}\right\}$, $\alpha \leftarrow \frac{1}{3\tau_1 L}$; ◇ *parameters*
$y_0 = z_0 = \widetilde{x}^0 \leftarrow x_0$; ◇ *initial vectors*
**for** $s \leftarrow 0$ **to** $S-1$ **do**
  $\mu^s \leftarrow \nabla f(\widetilde{x}^s)$; ◇ *compute the full gradient once every $m$ iterations*
  **for** $j \leftarrow 0$ **to** $m-1$ **do**
    $k \leftarrow (sm) + j$;
    $x_{k+1} \leftarrow \tau_1 z_k + \tau_2 \widetilde{x}^s + (1 - \tau_1 - \tau_2) y_k$;
    $\widetilde{\nabla}_{k+1} \leftarrow \mu^s + \nabla f_i(x_{k+1}) - \nabla f_i(\widetilde{x}^s)$ where $i$ is random from $\{1, 2, \dots, n\}$;
    $z_{k+1} = \arg\min_z \left\{\frac{1}{2\alpha}\|z - z_k\|^2 + \langle \widetilde{\nabla}_{k+1}, z\rangle + \psi(z)\right\}$;
    Option I: $y_{k+1} \leftarrow \arg\min_y \left\{\frac{3L}{2}\|y - x_{k+1}\|^2 + \langle \widetilde{\nabla}_{k+1}, y\rangle + \psi(y)\right\}$;
    Option II: $y_{k+1} \leftarrow x_{k+1} + \tau_1(z_{k+1} - z_k)$ ◇ *we analyze only I but II also works*
  **end for**
  $\widetilde{x}^{s+1} \leftarrow \left(\sum_{j=0}^{m-1}(1+\alpha\sigma)^j\right)^{-1} \cdot \left(\sum_{j=0}^{m-1}(1+\alpha\sigma)^j \cdot y_{sm+j+1}\right)$; ◇ *compute snapshot $\widetilde{x}$*
**end for**
**return** $\widetilde{x}^S$.

choices of gradient estimators. We choose $\widetilde{x}$ to be a weighted average, rather than the last or the uniform average, because it yields the tightest possible result.[9]

- *$\tau_1$ and $\alpha$ are standard parameters already present in Nesterov's full-gradient method [12].*

  We choose $\alpha = 1/3\tau_1 L$ to present the simplest proof, and recall it was $\alpha = 1/\tau_1 L$ in the original Nesterov's full-gradient method. (Any $\alpha$ that is constant factor smaller than $1/\tau_1 L$ works in theory, and we use 1/3 to provide the simplest proof.) In practice, like other accelerated methods, it suffices to fix $\alpha = 1/3\tau_1 L$ and only tune $\tau_1$ and thus $\tau_1$ is viewed as the learning rate.

- *The parameter $\tau_2$ is our novel weight for the Katyusha momentum. Any constant in $(0,1)$ works for $\tau_2$, and we simply choose $\tau_2 = 1/2$ for our theoretical and experimental results.*

We state our main theorem for `Katyusha` as follows:

**Theorem 2.1.** *If each $f_i(x)$ is convex, $L$-smooth, and $\psi(x)$ is $\sigma$-strongly convex in Problem (1.1), then* `Katyusha`$(x_0, S, \sigma, L)$ *satisfies*

$$\mathbb{E}\big[F(\widetilde{x}^S)\big] - F(x^*) \le \begin{cases} O\Big(\big(1 + \sqrt{\sigma/(3Lm)}\big)^{-Sm}\Big) \cdot \big(F(x_0) - F(x^*)\big), & \text{if } \frac{m\sigma}{L} \le \frac{3}{4}; \\ O\big(1.5^{-S}\big) \cdot \big(F(x_0) - F(x^*)\big), & \text{if } \frac{m\sigma}{L} > \frac{3}{4}. \end{cases}$$

*In other words, choosing $m = \Theta(n)$,* `Katyusha` *achieves an $\varepsilon$-additive error (i.e., $\mathbb{E}\big[F(\widetilde{x}^S)\big] - F(x^*) \le \varepsilon$) using at most $O\Big(\big(n + \sqrt{nL/\sigma}\big) \cdot \log \frac{F(x_0) - F(x^*)}{\varepsilon}\Big)$ iterations.*[10]

The proof of Theorem 2.1 is included in Section 2.1 and 2.2. As discussed in Section 1.1, the main idea behind our theorem is the negative momentum that helps reduce the error occurred from the stochastic gradient estimator.

*Remark* 2.2. Because $m = 2n$, each iteration of `Katyusha` computes only 1.5 stochastic gradients

[9] If one uses the uniform average, in theory, the algorithm needs to restart every a number of epochs (that is, by resetting $k = 0$, $s = 0$, and $x_0 = y_0 = z_0$); we refrain from doing so because we wish to provide a simple and direct algorithm. We can also use the last iterate, then the total complexity loses a factor $\log(L/\sigma)$. In practice, it was reported that even for SVRG, choosing average works better than choosing the last iterate [13].

$\nabla f_i(\cdot)$ in the amortized sense, the same as non-accelerated methods such as SVRG [29].[11] Therefore, the per-iteration cost of `Katyusha` is dominated by the computation of $\nabla f_i(\cdot)$, the proximal update in Line 10 of Algorithm 1, plus an overhead $O(d)$. If $\nabla f_i(\cdot)$ has at most $d' \leq d$ non-zero entries, this overhead $O(d)$ is improvable to $O(d')$ using a sparse implementation of `Katyusha`.[12]

For ERM problems defined in Problem (1.3), the amortized per-iteration complexity of `Katyusha` is $O(d')$ where $d'$ is the sparsity of feature vectors, the same as the per-iteration complexity of SGD.

## 2.1 One-Iteration Analysis

In this subsection, we first analyze the behavior of `Katyusha` in a single iteration (i.e., for a fixed $k$). We view $y_k, z_k$ and $x_{k+1}$ as fixed in this section so the only randomness comes from the choice of $i$ in iteration $k$. We abbreviate in this subsection by $\widetilde{x} = \widetilde{x}^s$ where $s$ is the epoch that iteration $k$ belongs to, and denote by $\sigma_{k+1}^2 \stackrel{\text{def}}{=} \|\nabla f(x_{k+1}) - \widetilde{\nabla}_{k+1}\|^2$ so $\mathbb{E}[\sigma_{k+1}^2]$ is the variance of the gradient estimator $\widetilde{\nabla}_{k+1}$ in this iteration.

Our first lemma lower bounds the expected objective decrease $F(x_{k+1}) - \mathbb{E}[F(y_{k+1})]$. Our $\mathsf{Prog}(x_{k+1})$ defined below is a non-negative, classical quantity that would be a lower bound on the amount of objective decrease if $\widetilde{\nabla}_{k+1}$ were equal to $\nabla f(x_{k+1})$ [12]. However, since the variance $\sigma_{k+1}^2$ is non-zero, this lower bound must be compensated by a negative term that depends on $\mathbb{E}[\sigma_{k+1}^2]$.

**Lemma 2.3** (proximal gradient descent). *If*

$$y_{k+1} = \arg\min_y \Big\{\frac{3L}{2}\|y - x_{k+1}\|^2 + \langle \widetilde{\nabla}_{k+1}, y - x_{k+1}\rangle + \psi(y) - \psi(x_{k+1})\Big\} \ , \qquad \textit{and}$$

$$\mathsf{Prog}(x_{k+1}) \stackrel{\text{def}}{=} -\min_y \Big\{\frac{3L}{2}\|y - x_{k+1}\|^2 + \langle \widetilde{\nabla}_{k+1}, y - x_{k+1}\rangle + \psi(y) - \psi(x_{k+1})\Big\} \geq 0 \ ,$$

*we have (where the expectation is only over the randomness of* $\widetilde{\nabla}_{k+1}$*)*

$$F(x_{k+1}) - \mathbb{E}\big[F(y_{k+1})\big] \geq \mathbb{E}\big[\mathsf{Prog}(x_{k+1})\big] - \frac{1}{4L}\mathbb{E}\big[\sigma_{k+1}^2\big] \ .$$

*Proof.*

$$\begin{aligned}
\mathsf{Prog}(x_{k+1}) &= -\min_y \Big\{\frac{3L}{2}\|y - x_{k+1}\|^2 + \langle \widetilde{\nabla}_{k+1}, y - x_{k+1}\rangle + \psi(y) - \psi(x_{k+1})\Big\} \\
&\overset{①}{=} -\Big(\frac{3L}{2}\|y_{k+1} - x_{k+1}\|^2 + \langle \widetilde{\nabla}_{k+1}, y_{k+1} - x_{k+1}\rangle + \psi(y_{k+1}) - \psi(x_{k+1})\Big) \\
&= -\Big(\frac{L}{2}\|y_{k+1} - x_{k+1}\|^2 + \langle \nabla f(x_{k+1}), y_{k+1} - x_{k+1}\rangle + \psi(y_{k+1}) - \psi(x_{k+1})\Big) \\
&\quad + \Big(\langle \nabla f(x_{k+1}) - \widetilde{\nabla}_{k+1}, y_{k+1} - x_{k+1}\rangle - L\|y_{k+1} - x_{k+1}\|^2\Big) \\
&\overset{②}{\leq} -\Big(f(y_{k+1}) - f(x_{k+1}) + \psi(y_{k+1}) - \psi(x_{k+1})\Big) + \frac{1}{4L}\|\nabla f(x_{k+1}) - \widetilde{\nabla}_{k+1}\|^2 \ .
\end{aligned}$$

Above, ① is by the definition of $y_{k+1}$, and ② uses the smoothness of function $f(\cdot)$, as well as Young's inequality $\langle a, b\rangle - \frac{1}{2}\|b\|^2 \leq \frac{1}{2}\|a\|^2$. Taking expectation on both sides we arrive at the desired result. □

[11] The claim "SVRG or `Katyusha` computes 1.5 stochastic gradients" requires one to store $\nabla_i f(\widetilde{x})$ in the memory for each $i \in [n]$, and this costs $O(dn)$ space in the most general setting. If one does not store $\nabla_i f(\widetilde{x})$ in the memory, then each iteration of SVRG or `Katyusha` computes 2.5 stochastic gradients for the choice $m = 2n$.

[12] This requires to defer a coordinate update to the moment it is accessed. Update deferral is a standard technique used in sparse implementations of all stochastic gradient methods, including SVRG, SAGA, APCG [22, 29, 37].

The following lemma provides a novel upper bound on the expected variance of the gradient estimator. Note that all known variance reduction analysis for convex optimization, in one way or another, upper bounds this variance essentially by $4L \cdot (f(\widetilde{x}) - f(x^*))$, the objective distance to the minimizer (c.f. [22, 29]). The recent result of Allen-Zhu and Hazan [6] upper bounds it by the point distance $\|x_{k+1} - \widetilde{x}\|^2$ for non-convex objectives, which is tighter if $\widetilde{x}$ is close to $x_{k+1}$ but unfortunately not enough for the purpose of this paper.

In this paper, we upper bound it by the tightest possible quantity which is essentially $2L \cdot \big(f(\widetilde{x}) - f(x_{k+1})\big) \ll 4L \cdot \big(f(\widetilde{x}) - f(x^*)\big)$. Unfortunately, this upper bound needs to be compensated by an additional term $\langle \nabla f(x_{k+1}), \widetilde{x} - x_{k+1}\rangle$, which could be positive but we shall cancel it using the introduced Katyusha momentum.

**Lemma 2.4** (variance upper bound)**.**

$$\mathbb{E}\big[\|\widetilde{\nabla}_{k+1} - \nabla f(x_{k+1})\|^2\big] \leq 2L \cdot \big(f(\widetilde{x}) - f(x_{k+1}) - \langle \nabla f(x_{k+1}), \widetilde{x} - x_{k+1}\rangle\big) \ .$$

*Proof.* Each $f_i(x)$, being convex and $L$-smooth, implies the following inequality which is classical in convex optimization and can be found for instance in Theorem 2.1.5 of the textbook of Nesterov [43].

$$\|\nabla f_i(x_{k+1}) - \nabla f_i(\widetilde{x})\|^2 \leq 2L \cdot \big(f_i(\widetilde{x}) - f_i(x_{k+1}) - \langle \nabla f_i(x_{k+1}), \widetilde{x} - x_{k+1}\rangle\big)$$

Therefore, taking expectation over the random choice of $i$, we have

$$\begin{aligned}
\mathbb{E}\big[\|\widetilde{\nabla}_{k+1} - \nabla f(x_{k+1})\|^2\big] &= \mathbb{E}\Big[\big\|\big(\nabla f_i(x_{k+1}) - \nabla f_i(\widetilde{x})\big) - \big(\nabla f(x_{k+1}) - \nabla f(\widetilde{x})\big)\big\|^2\Big] \\
&\overset{①}{\leq} \mathbb{E}\big[\big\|\nabla f_i(x_{k+1}) - \nabla f_i(\widetilde{x})\big\|^2\big] \\
&\overset{②}{\leq} 2L \cdot \mathbb{E}\big[f_i(\widetilde{x}) - f_i(x_{k+1}) - \langle \nabla f_i(x_{k+1}), \widetilde{x} - x_{k+1}\rangle\big] \\
&= 2L \cdot \big(f(\widetilde{x}) - f(x_{k+1}) - \langle \nabla f(x_{k+1}), \widetilde{x} - x_{k+1}\rangle\big) \ .
\end{aligned}$$

Above, ① is because for any random vector $\zeta \in \mathbb{R}^d$, it holds that $\mathbb{E}\|\zeta - \mathbb{E}\zeta\|^2 = \mathbb{E}\|\zeta\|^2 - \|\mathbb{E}\zeta\|^2$; ② follows from the first inequality in this proof. □

The next lemma is a classical one for proximal mirror descent.

**Lemma 2.5** (proximal mirror descent)**.** *Suppose $\psi(\cdot)$ is $\sigma$-SC. Then, fixing $\widetilde{\nabla}_{k+1}$ and letting*

$$z_{k+1} = \arg\min_{z} \Big\{\frac{1}{2}\|z - z_k\|^2 + \alpha\langle \widetilde{\nabla}_{k+1}, z - z_k\rangle + \alpha\psi(z) - \alpha\psi(z_k)\Big\} \ ,$$

*it satisfies for all $u \in \mathbb{R}^d$,*

$$\alpha\langle \widetilde{\nabla}_{k+1}, z_{k+1} - u\rangle + \alpha\psi(z_{k+1}) - \alpha\psi(u) \leq -\frac{1}{2}\|z_k - z_{k+1}\|^2 + \frac{1}{2}\|z_k - u\|^2 - \frac{1+\alpha\sigma}{2}\|z_{k+1} - u\|^2 \ .$$

*Proof.* By the minimality definition of $z_{k+1}$, we have that

$$z_{k+1} - z_k + \alpha\widetilde{\nabla}_{k+1} + \alpha g = 0$$

where $g$ is *some* subgradient of $\psi(z)$ at point $z = z_{k+1}$. This implies that for every $u$ it satisfies

$$0 = \big\langle z_{k+1} - z_k + \alpha\widetilde{\nabla}_{k+1} + \alpha g, z_{k+1} - u\big\rangle \ .$$

At this point, using the equality $\langle z_{k+1} - z_k, z_{k+1} - u\rangle = \frac{1}{2}\|z_k - z_{k+1}\|^2 - \frac{1}{2}\|z_k - u\|^2 + \frac{1}{2}\|z_{k+1} - u\|^2$, as well as the inequality $\langle g, z_{k+1} - u\rangle \geq \psi(z_{k+1}) - \psi(u) + \frac{\sigma}{2}\|z_{k+1} - u\|^2$ which comes from the strong

convexity of $\psi(\cdot)$, we can write

$$\begin{aligned}
&\alpha\langle\widetilde{\nabla}_{k+1}, z_{k+1}-u\rangle + \alpha\psi(z_{k+1}) - \alpha\psi(u) \\
&= -\langle z_{k+1}-z_k, z_{k+1}-u\rangle - \langle \alpha g, z_{k+1}-u\rangle + \alpha\psi(z_{k+1}) - \alpha\psi(u) \\
&\le -\frac{1}{2}\|z_k - z_{k+1}\|^2 + \frac{1}{2}\|z_k-u\|^2 - \frac{1+\alpha\sigma}{2}\|z_{k+1}-u\|^2 \enspace . \qquad \square
\end{aligned}$$

The following lemma combines Lemma 2.3, Lemma 2.4 and Lemma 2.5 all together, using the special choice of $x_{k+1}$ which is a convex combination of $y_k, z_k$ and $\widetilde{x}$:

**Lemma 2.6** (coupling step 1). *If $x_{k+1} = \tau_1 z_k + \tau_2\widetilde{x} + (1-\tau_1-\tau_2)y_k$, where $\tau_1 \le \frac{1}{3\alpha L}$ and $\tau_2 = \frac{1}{2}$,*

$$\begin{aligned}
&\alpha\langle\nabla f(x_{k+1}), z_k - u\rangle - \alpha\psi(u) \\
&\le \frac{\alpha}{\tau_1}\Big(F(x_{k+1}) - \mathbb{E}\big[F(y_{k+1})\big] + \tau_2 F(\widetilde{x}) - \tau_2 f(x_{k+1}) - \tau_2\langle\nabla f(x_{k+1}), \widetilde{x} - x_{k+1}\rangle\Big) \\
&\quad + \frac{1}{2}\|z_k-u\|^2 - \frac{1+\alpha\sigma}{2}\mathbb{E}\big[\|z_{k+1}-u\|^2\big] + \frac{\alpha(1-\tau_1-\tau_2)}{\tau_1}\psi(y_k) - \frac{\alpha}{\tau_1}\psi(x_{k+1}) \enspace .
\end{aligned}$$

*Proof.* We first apply Lemma 2.5 and get

$$\begin{aligned}
&\alpha\langle\widetilde{\nabla}_{k+1}, z_k - u\rangle + \alpha\psi(z_{k+1}) - \alpha\psi(u) \\
&= \alpha\langle\widetilde{\nabla}_{k+1}, z_k - z_{k+1}\rangle + \alpha\langle\widetilde{\nabla}_{k+1}, z_{k+1} - u\rangle + \alpha\psi(z_{k+1}) - \alpha\psi(u) \\
&\le \alpha\langle\widetilde{\nabla}_{k+1}, z_k - z_{k+1}\rangle - \frac{1}{2}\|z_k - z_{k+1}\|^2 + \frac{1}{2}\|z_k-u\|^2 - \frac{1+\alpha\sigma}{2}\|z_{k+1}-u\|^2 \enspace . \qquad (2.1)
\end{aligned}$$

By defining $v \stackrel{\text{def}}{=} \tau_1 z_{k+1} + \tau_2\widetilde{x} + (1-\tau_1-\tau_2)y_k$, we have $x_{k+1} - v = \tau_1(z_k - z_{k+1})$ and therefore

$$\begin{aligned}
&\mathbb{E}\Big[\alpha\langle\widetilde{\nabla}_{k+1}, z_k - z_{k+1}\rangle - \frac{1}{2}\|z_k - z_{k+1}\|^2\Big] = \mathbb{E}\Big[\frac{\alpha}{\tau_1}\langle\widetilde{\nabla}_{k+1}, x_{k+1} - v\rangle - \frac{1}{2\tau_1^2}\|x_{k+1}-v\|^2\Big] \\
&= \mathbb{E}\Big[\frac{\alpha}{\tau_1}\Big(\langle\widetilde{\nabla}_{k+1}, x_{k+1}-v\rangle - \frac{1}{2\alpha\tau_1}\|x_{k+1}-v\|^2 - \psi(v) + \psi(x_{k+1})\Big) + \frac{\alpha}{\tau_1}\Big(\psi(v) - \psi(x_{k+1})\Big)\Big] \\
&\overset{①}{\le} \mathbb{E}\Big[\frac{\alpha}{\tau_1}\Big(\langle\widetilde{\nabla}_{k+1}, x_{k+1}-v\rangle - \frac{3L}{2}\|x_{k+1}-v\|^2 - \psi(v) + \psi(x_{k+1})\Big) + \frac{\alpha}{\tau_1}\Big(\psi(v) - \psi(x_{k+1})\Big)\Big] \\
&\overset{②}{\le} \mathbb{E}\Big[\frac{\alpha}{\tau_1}\Big(F(x_{k+1}) - F(y_{k+1}) + \frac{1}{4L}\sigma_{k+1}^2\Big) + \frac{\alpha}{\tau_1}\Big(\psi(v) - \psi(x_{k+1})\Big)\Big] \\
&\overset{③}{\le} \mathbb{E}\Big[\frac{\alpha}{\tau_1}\Big(F(x_{k+1}) - F(y_{k+1}) + \frac{1}{2}\big(f(\widetilde{x}) - f(x_{k+1}) - \langle\nabla f(x_{k+1}), \widetilde{x} - x_{k+1}\rangle\big)\Big) \\
&\quad + \frac{\alpha}{\tau_1}\Big(\tau_1\psi(z_{k+1}) + \tau_2\psi(\widetilde{x}) + (1-\tau_1-\tau_2)\psi(y_k) - \psi(x_{k+1})\Big)\Big] \enspace . \qquad (2.2)
\end{aligned}$$

Above, ① uses our choice $\tau_1 \le \frac{1}{3\alpha L}$, ② uses Lemma 2.3, ③ uses Lemma 2.4 together with the convexity of $\psi(\cdot)$ and the definition of $v$. Finally, noticing that $\mathbb{E}[\langle\widetilde{\nabla}_{k+1}, z_k - u\rangle] = \langle\nabla f(x_{k+1}), z_k - u\rangle$ and $\tau_2 = \frac{1}{2}$, we obtain the desired inequality by combining (2.1) and (2.2). $\square$

The next lemma simplifies the left hand side of Lemma 2.6 using the convexity of $f(\cdot)$, and gives an inequality that relates the objective-distance-to-minimizer quantities $F(y_k) - F(x^*)$, $F(y_{k+1}) - F(x^*)$, and $F(\widetilde{x}) - F(x^*)$ to the point-distance-to-minimizer quantities $\|z_k - x^*\|^2$ and $\|z_{k+1} - x^*\|^2$.

**Lemma 2.7** (coupling step 2)**.** *Under the same choices of $\tau_1, \tau_2$ as in Lemma 2.6, we have*

$$0 \le \frac{\alpha(1-\tau_1-\tau_2)}{\tau_1}(F(y_k)-F(x^*)) - \frac{\alpha}{\tau_1}\big(\mathbb{E}\big[F(y_{k+1})\big]-F(x^*)\big) + \frac{\alpha\tau_2}{\tau_1}\big(F(\widetilde{x})-F(x^*)\big)$$
$$+\frac{1}{2}\|z_k-x^*\|^2 - \frac{1+\alpha\sigma}{2}\mathbb{E}\big[\|z_{k+1}-x^*\|^2\big] \ .$$

*Proof.* We first compute that

$$\begin{aligned}
&\alpha\big(f(x_{k+1})-f(u)\big) \overset{①}{\le} \alpha\langle\nabla f(x_{k+1}), x_{k+1}-u\rangle \\
&= \alpha\langle\nabla f(x_{k+1}), x_{k+1}-z_k\rangle + \alpha\langle\nabla f(x_{k+1}), z_k-u\rangle \\
&\overset{②}{=} \frac{\alpha\tau_2}{\tau_1}\langle\nabla f(x_{k+1}), \widetilde{x}-x_{k+1}\rangle + \frac{\alpha(1-\tau_1-\tau_2)}{\tau_1}\langle\nabla f(x_{k+1}), y_k-x_{k+1}\rangle + \alpha\langle\nabla f(x_{k+1}), z_k-u\rangle \\
&\overset{③}{\le} \frac{\alpha\tau_2}{\tau_1}\langle\nabla f(x_{k+1}), \widetilde{x}-x_{k+1}\rangle + \frac{\alpha(1-\tau_1-\tau_2)}{\tau_1}(f(y_k)-f(x_{k+1})) + \alpha\langle\nabla f(x_{k+1}), z_k-u\rangle \ .
\end{aligned}$$

Above, ① uses the convexity of $f(\cdot)$, ② uses the choice that $x_{k+1} = \tau_1 z_k + \tau_2\widetilde{x} + (1-\tau_1-\tau_2)y_k$, and ③ uses the convexity of $f(\cdot)$ again. By applying Lemma 2.6 to the above inequality, we have

$$\alpha\big(f(x_{k+1})-F(u)\big) \le \frac{\alpha(1-\tau_1-\tau_2)}{\tau_1}(F(y_k)-f(x_{k+1}))$$
$$+\frac{\alpha}{\tau_1}\Big(F(x_{k+1})-\mathbb{E}\big[F(y_{k+1})\big]+\tau_2F(\widetilde{x})-\tau_2f(x_{k+1})\Big)+\frac{1}{2}\|z_k-u\|^2-\frac{1+\alpha\sigma}{2}\mathbb{E}\big[\|z_{k+1}-u\|^2\big]-\frac{\alpha}{\tau_1}\psi(x_{k+1})$$

which implies

$$\alpha\big(F(x_{k+1})-F(u)\big) \le \frac{\alpha(1-\tau_1-\tau_2)}{\tau_1}(F(y_k)-F(x_{k+1}))$$
$$+\frac{\alpha}{\tau_1}\Big(F(x_{k+1})-\mathbb{E}\big[F(y_{k+1})\big]+\tau_2F(\widetilde{x})-\tau_2F(x_{k+1})\Big)+\frac{1}{2}\|z_k-u\|^2-\frac{1+\alpha\sigma}{2}\mathbb{E}\big[\|z_{k+1}-u\|^2\big] \ .$$

After rearranging and setting $u = x^*$, the above inequality yields

$$0 \le \frac{\alpha(1-\tau_1-\tau_2)}{\tau_1}(F(y_k)-F(x^*)) - \frac{\alpha}{\tau_1}\big(\mathbb{E}\big[F(y_{k+1})-F(x^*)\big]\big) + \frac{\alpha\tau_2}{\tau_1}\big(F(\widetilde{x})-F(x^*)\big)$$
$$+\frac{1}{2}\|z_k-x^*\|^2 - \frac{1+\alpha\sigma}{2}\mathbb{E}\big[\|z_{k+1}-x^*\|^2\big] \ . \qquad \square$$

## 2.2 Proof of Theorem 2.1

We are now ready to combine the analyses across iterations, and derive our final Theorem 2.1. Our proof next requires a careful telescoping of Lemma 2.7 together with our specific parameter choices.

*Proof of Theorem 2.1.* Define $D_k \overset{\text{def}}{=} F(y_k)-F(x^*)$, $\widetilde{D}^s \overset{\text{def}}{=} F(\widetilde{x}^s)-F(x^*)$, and rewrite Lemma 2.7:

$$0 \le \frac{(1-\tau_1-\tau_2)}{\tau_1}D_k - \frac{1}{\tau_1}D_{k+1} + \frac{\tau_2}{\tau_1}\mathbb{E}\big[\widetilde{D}^s\big] + \frac{1}{2\alpha}\|z_k-x^*\|^2 - \frac{1+\alpha\sigma}{2\alpha}\mathbb{E}\big[\|z_{k+1}-x^*\|^2\big] \ .$$

At this point, let us define $\theta = 1+\alpha\sigma$ and multiply the above inequality by $\theta^j$ for each $k = sm+j$. Then, we sum up the resulting $m$ inequalities for all $j = 0, 1, \dots, m-1$:

$$0 \le \mathbb{E}\Big[\frac{(1-\tau_1-\tau_2)}{\tau_1}\sum_{j=0}^{m-1} D_{sm+j}\cdot\theta^j - \frac{1}{\tau_1}\sum_{j=0}^{m-1} D_{sm+j+1}\cdot\theta^j\Big] + \frac{\tau_2}{\tau_1}\widetilde{D}^s\cdot\sum_{j=0}^{m-1}\theta^j$$
$$+\frac{1}{2\alpha}\|z_{sm}-x^*\|^2 - \frac{\theta^m}{2\alpha}\big[\|z_{(s+1)m}-x^*\|^2\big] \ .$$

Note that in the above inequality we have assumed all the randomness in the first $s-1$ epochs are fixed and the only source of randomness comes from epoch $s$. We can rearrange the terms in the above inequality and get

$$\mathbb{E}\Big[\frac{\tau_1+\tau_2-(1-1/\theta)}{\tau_1}\sum_{j=1}^{m} D_{sm+j}\cdot\theta^j\Big] \leq \frac{(1-\tau_1-\tau_2)}{\tau_1}\Big(D_{sm}-\theta^m\mathbb{E}\big[D_{(s+1)m}\big]\Big)$$
$$+\frac{\tau_2}{\tau_1}\widetilde{D}^s\cdot\sum_{j=0}^{m-1}\theta^j+\frac{1}{2\alpha}\|z_{sm}-x^*\|^2-\frac{\theta^m}{2\alpha}\mathbb{E}\big[\|z_{(s+1)m}-x^*\|^2\big] \enspace .$$

Using the special choice that $\widetilde{x}^{s+1}=\big(\sum_{j=0}^{m-1}\theta^j\big)^{-1}\cdot\sum_{j=0}^{m-1}y_{sm+j+1}\cdot\theta^j$ and the convexity of $F(\cdot)$, we derive that $\widetilde{D}^{s+1}\leq\big(\sum_{j=0}^{m-1}\theta^j\big)^{-1}\cdot\sum_{j=0}^{m-1}D_{sm+j+1}\cdot\theta^j$. Substituting this into the above inequality, we get

$$\frac{\tau_1+\tau_2-(1-1/\theta)}{\tau_1}\theta\mathbb{E}\big[\widetilde{D}^{s+1}\big]\cdot\sum_{j=0}^{m-1}\theta^j \leq \frac{(1-\tau_1-\tau_2)}{\tau_1}\Big(D_{sm}-\theta^m\mathbb{E}\big[D_{(s+1)m}\big]\Big)$$
$$+\frac{\tau_2}{\tau_1}\widetilde{D}^s\cdot\sum_{j=0}^{m-1}\theta^j+\frac{1}{2\alpha}\|z_{sm}-x^*\|^2-\frac{\theta^m}{2\alpha}\mathbb{E}\big[\|z_{(s+1)m}-x^*\|^2\big] \enspace . \quad (2.3)$$

We consider two cases next.

**Case 1.** Suppose $\frac{m\sigma}{L}\leq\frac{3}{4}$. In this case, we choose $\alpha=\frac{1}{\sqrt{3m\sigma L}}$ and $\tau_1=\frac{1}{3\alpha L}=m\alpha\sigma=\frac{\sqrt{m\sigma}}{\sqrt{3L}}\in[0,\frac{1}{2}]$ for `Katyusha`. It implies $\alpha\sigma\leq 1/2m$ and therefore the following inequality holds:

$$\tau_2(\theta^{m-1}-1)+(1-1/\theta)=\frac{1}{2}((1+\alpha\sigma)^{m-1}-1)+(1-\frac{1}{1+\alpha\sigma})\leq(m-1)\alpha\sigma+\alpha\sigma=m\alpha\sigma=\tau_1 \enspace .$$

In other words, we have $\tau_1+\tau_2-(1-1/\theta)\geq\tau_2\theta^{m-1}$ and thus (2.3) implies that

$$\mathbb{E}\Big[\frac{\tau_2}{\tau_1}\widetilde{D}^{s+1}\cdot\sum_{j=0}^{m-1}\theta^j+\frac{1-\tau_1-\tau_2}{\tau_1}D_{(s+1)m}+\frac{1}{2\alpha}\|z_{(s+1)m}-x^*\|^2\Big]$$
$$\leq\theta^{-m}\cdot\Big(\frac{\tau_2}{\tau_1}\widetilde{D}^s\cdot\sum_{j=0}^{m-1}\theta^j+\frac{1-\tau_1-\tau_2}{\tau_1}D_{sm}+\frac{1}{2\alpha}\|z_{sm}-x^*\|^2\Big) \enspace .$$

If we telescope the above inequality over all epochs $s=0,1,\dots,S-1$, we obtain

$$\mathbb{E}\big[F(\widetilde{x}^S)-F(x^*)\big]=\mathbb{E}\big[\widetilde{D}^S\big]\overset{①}{\leq}\theta^{-Sm}\cdot O\Big(\widetilde{D}^0+D_0+\frac{\tau_1}{\alpha m}\|x_0-x^*\|^2\Big)$$
$$\overset{②}{\leq}\theta^{-Sm}\cdot O\Big(1+\frac{\tau_1}{\alpha m\sigma}\Big)\cdot(F(x_0)-F(x^*))$$
$$\overset{③}{=}O((1+\alpha\sigma)^{-Sm})\cdot\big(F(x_0)-F(x^*)\big) \enspace . \quad (2.4)$$

Above, ① uses the fact that $\sum_{j=0}^{m-1}\theta^j\geq m$ and $\tau_2=\frac{1}{2}$; ② uses the strong convexity of $F(\cdot)$ which implies $F(x_0)-F(x^*)\geq\frac{\sigma}{2}\|x_0-x^*\|^2$; and ③ uses our choice of $\tau_1$.

**Case 2.** Suppose $\frac{m\sigma}{L}>\frac{3}{4}$. In this case, we choose $\tau_1=\frac{1}{2}$ and $\alpha=\frac{1}{3\tau_1 L}=\frac{2}{3L}$ as in `Katyusha`. Our parameter choices help us simplify (2.3) as (noting $(\tau_1+\tau_2-(1-1/\theta))\theta=1$)

$$2\mathbb{E}\big[\widetilde{D}^{s+1}\big]\cdot\sum_{j=0}^{m-1}\theta^j\leq\widetilde{D}^s\cdot\sum_{j=0}^{m-1}\theta^j+\frac{1}{2\alpha}\|z_{sm}-x^*\|^2-\frac{\theta^m}{2\alpha}\mathbb{E}\big[\|z_{(s+1)m}-x^*\|^2\big] \enspace .$$

Since $\theta^m = (1+\alpha\sigma)^m \geq 1 + \alpha\sigma m = 1 + \frac{2\sigma m}{3L} \geq \frac{3}{2}$, the above inequality implies

$$\frac{3}{2}\mathbb{E}\big[\widetilde{D}^{s+1}\big] \cdot \sum_{j=0}^{m-1} \theta^j + \frac{9L}{8}\mathbb{E}\big[\|z_{(s+1)m} - x^*\|^2\big] \leq \widetilde{D}^s \cdot \sum_{j=0}^{m-1} \theta^j + \frac{3L}{4}\|z_{sm} - x^*\|^2 \enspace .$$

If we telescope this inequality over all the epochs $s = 0, 1, \dots, S-1$, we immediately have

$$\mathbb{E}\Big[\widetilde{D}^S \cdot \sum_{j=0}^{m-1} \theta^j + \frac{3L}{4}\|z_{Sm} - x^*\|^2\Big] \leq \big(\frac{2}{3}\big)^S \cdot \Big(\widetilde{D}^0 \cdot \sum_{j=0}^{m-1} \theta^j + \frac{3L}{4}\|z_0 - x^*\|^2\Big) \enspace .$$

Finally, since $\sum_{j=0}^{m-1} \theta^j \geq m$ and $\frac{\sigma}{2}\|z_0 - x^*\|^2 \leq F(x_0) - F(x^*)$ owing to the strong convexity of $F(\cdot)$, we conclude that

$$\mathbb{E}\big[F(\widetilde{x}^S) - F(x^*)\big] \leq O\big(1.5^{-S}\big) \cdot \big(F(x_0) - F(x^*)\big) \enspace . \tag{2.5}$$

Combining (2.4) and (2.5) we finish the proof of Theorem 2.1. □

# 3 Corollaries on Non-Smooth or Non-SC Problems

In this section we apply reductions to translate our Theorem 2.1 into optimal algorithms also for non-strongly convex objectives and/or non-smooth objectives.

To begin with, recall the following definition of the HOOD property:

**Definition 3.1** ([6]). *An algorithm solving the strongly convex case of Problem (1.1) satisfies the* homogenous objective decrease (HOOD) *property with $T(L,\sigma)$, if for every starting point $x_0$, it produces an output $x'$ satisfying $\mathbb{E}\big[F(x')\big] - F(x^*) \leq \frac{F(x_0) - F(x^*)}{4}$ in at most $T(L,\sigma)$ stochastic gradient iterations.*

Theorem 2.1 shows that `Katyusha` satisfies the HOOD property:

**Corollary 3.2.** `Katyusha` *satisfies the HOOD property with $T(L,\sigma) = O\big(n + \frac{\sqrt{nL}}{\sqrt{\sigma}}\big)$.*

*Remark* 3.3. Existing accelerated stochastic methods before this work (even for simpler Problem (1.3)) either do not satisfy HOOD or satisfy HOOD with an additional factor $\log(L/\sigma)$ in the number of iterations.

Allen-Zhu and Hazan [6] designed three reductions algorithms to convert an algorithm satisfying the HOOD property to solve the following three cases:

**Theorem 3.4.** *Given algorithm $\mathcal{A}$ satisfying HOOD with $T(L,\sigma)$ and a starting vector $x_0$.*

- NonSC+Smooth. *For Problem (1.1) where $f(\cdot)$ is $L$-smooth,* AdaptReg($\mathcal{A}$) *outputs $x$ satisfying $\mathbb{E}\big[F(x)\big] - F(x^*) \leq O(\varepsilon)$ in $T$ stochastic gradient iterations where*

$$T = \sum_{s=0}^{S-1} T\Big(L, \frac{\sigma_0}{2^s}\Big) \text{ where } \sigma_0 = \frac{F(x_0) - F(x^*)}{\|x_0 - x^*\|^2} \text{ and } S = \log_2 \frac{F(x_0) - F(x^*)}{\varepsilon}.$$

- SC+NonSmooth. *For Problem (1.3) where $\psi(\cdot)$ is $\sigma$-SC and each $f_i(\cdot)$ is $\sqrt{G}$-Lipschitz continuous,* AdaptSmooth($\mathcal{A}$) *outputs $x$ satisfying $\mathbb{E}\big[F(x)\big] - F(x^*) \leq O(\varepsilon)$ in*

$$T = \sum_{s=0}^{S-1} T\Big(\frac{2^s}{\lambda_0}, \sigma\Big) \text{ where } \lambda_0 = \frac{F(x_0) - F(x^*)}{G} \text{ and } S = \log_2 \frac{F(x_0) - F(x^*)}{\varepsilon}.$$

- NonSC+NonSmooth. *For Problem (1.3) where each $f_i(\cdot)$ is $\sqrt{G}$-Lipschitz continuous, then* $\mathsf{JointAdaptRegSmooth}(\mathcal{A})$ *outputs $x$ satisfying $\mathbb{E}\big[F(x)\big] - F(x^*) \le O(\varepsilon)$ in*

$$T = \sum_{s=0}^{S-1} T\Big(\frac{2^s}{\lambda_0}, \frac{\sigma_0}{2^s}\Big)$$

$$\text{where } \lambda_0 = \frac{F(x_0) - F(x^*)}{G},\ \sigma_0 = \frac{F(x_0) - F(x^*)}{\|x_0 - x^*\|^2} \text{ and } S = \log_2 \frac{F(x_0) - F(x^*)}{\|x_0 - x^*\|^2}.$$

Combining Corollary 3.2 with Theorem 3.4, we have the following corollaries:

**Corollary 3.5.** *If each $f_i(x)$ is convex, $L$-smooth and $\psi(\cdot)$ is not necessarily strongly convex in Problem (1.1), then by applying* AdaptReg *on* `Katyusha` *with a starting vector $x_0$, we obtain an output $x$ satisfying $\mathbb{E}[F(x)] - F(x^*) \le \varepsilon$ in*

$$T = O\Big(n \log \tfrac{F(x_0)-F(x^*)}{\varepsilon} + \tfrac{\sqrt{nL}\cdot\|x_0-x^*\|}{\sqrt{\varepsilon}}\Big) \propto \tfrac{1}{\sqrt{\varepsilon}} \text{ iterations. (Or equivalently } \varepsilon \propto \tfrac{1}{T^2}.)$$

In contrast, the best known convergence rate was $\varepsilon \propto \frac{\log^4 T}{T^2}$ or more precisely

$$\text{Catalyst:}\quad T = O\Big(\Big(n + \tfrac{\sqrt{nL}\cdot\|x_0-x^*\|}{\sqrt{\varepsilon}}\Big) \log \tfrac{F(x_0)-F(x^*)}{\varepsilon} \log \tfrac{L\|x_0-x^*\|^2}{\varepsilon}\Big) \propto \tfrac{\log^2(1/\varepsilon)}{\sqrt{\varepsilon}} \text{ iterations.}$$

**Corollary 3.6.** *If each $f_i(x)$ is $\sqrt{G}$-Lipschitz continuous and $\psi(x)$ is $\sigma$-SC in Problem (1.3), then by applying* AdaptSmooth *on* `Katyusha` *with a starting vector $x_0$, we obtain an output $x$ satisfying $\mathbb{E}[F(x)] - F(x^*) \le \varepsilon$ in*

$$T = O\Big(n \log \tfrac{F(x_0)-F(x^*)}{\varepsilon} + \tfrac{\sqrt{nG}}{\sqrt{\sigma\varepsilon}}\Big) \propto \tfrac{1}{\sqrt{\varepsilon}} \text{ iterations. (Or equivalently } \varepsilon \propto \tfrac{1}{T^2}.)$$

In contrast, the best known convergence rate was $\varepsilon \propto \frac{\log^2 T}{T^2}$, or more precisely

$$\text{APCG/SPDC:}\quad T = O\Big(\Big(n + \tfrac{\sqrt{nG}}{\sqrt{\sigma\varepsilon}}\Big) \log \tfrac{nG(F(x_0)-F(x^*))}{\sigma\varepsilon}\Big) \propto \tfrac{\log(1/\varepsilon)}{\sqrt{\varepsilon}} \text{ iterations.}$$

**Corollary 3.7.** *If each $f_i(x)$ is $\sqrt{G}$-Lipschitz continuous and $\psi(x)$ is not necessarily strongly convex in Problem (1.3), then by applying* JointAdaptRegSmooth *on* `Katyusha` *with a starting vector $x_0$, we obtain an output $x$ satisfying $\mathbb{E}[F(x)] - F(x^*) \le \varepsilon$ in*

$$T = O\Big(n \log \tfrac{F(x_0)-F(x^*)}{\varepsilon} + \tfrac{\sqrt{nG}\|x_0-x^*\|}{\varepsilon}\Big) \propto \tfrac{1}{\varepsilon} \text{ iterations. (Or equivalently } \varepsilon \propto \tfrac{1}{T}.)$$

In contrast, the best known convergence rate was $\varepsilon \propto \frac{\log T}{T}$, or more precisely

$$\text{APCG/SPDC:}\quad T = O\Big(\Big(n + \tfrac{\sqrt{nG}\|x_0-x^*\|}{\varepsilon}\Big) \log \tfrac{nG\|x_0-x^*\|^2(F(x_0)-F(x^*))}{\varepsilon^2}\Big) \propto \tfrac{\log(1/\varepsilon)}{\varepsilon} \text{ iterations.}$$

# 4 Katyusha in the Non-Strongly Convex Setting

Due to the increasing popularity of *non-strongly convex* minimization tasks (most notably $\ell_1$-regularized problems), researchers often make additional efforts to design separate methods for minimizing the non-strongly convex variant of Problem (1.1) that are *direct*, meaning without restarting and in particular without using any reductions such as Theorem 3.4 [13, 22].

In this section, we also develop our *direct and accelerated* method for the non-strongly convex variant of Problem (1.1). We call it $\texttt{Katyusha}^{\mathsf{ns}}$ and state it in Algorithm 2.

**Algorithm 2** $\texttt{Katyusha}^{\textsf{ns}}(x_0, S, \sigma, L)$

$m \leftarrow 2n$; ◇ *epoch length*
$\tau_2 \leftarrow \frac{1}{2}$;
$y_0 = z_0 = \widetilde{x}^0 \leftarrow x_0$; ◇ *initial vectors*
**for** $s \leftarrow 0$ **to** $S-1$ **do**
  $\tau_{1,s} \leftarrow \frac{2}{s+4}$, $\alpha_s \leftarrow \frac{1}{3\tau_{1,s}L}$ ◇ *different parameter choices comparing to* `Katyusha`
  $\mu^s \leftarrow \nabla f(\widetilde{x}^s)$; ◇ *compute the full gradient only once every m iterations*
  **for** $j \leftarrow 0$ **to** $m-1$ **do**
    $k \leftarrow (sm) + j$;
    $x_{k+1} \leftarrow \tau_{1,s} z_k + \tau_2 \widetilde{x}^s + (1 - \tau_{1,s} - \tau_2) y_k$;
    $\widetilde{\nabla}_{k+1} \leftarrow \mu^s + \nabla f_i(x_{k+1}) - \nabla f_i(\widetilde{x}^s)$ where $i$ is randomly chosen from $\{1, 2, \dots, n\}$;
    $z_{k+1} = \arg\min_z \big\{ \frac{1}{2\alpha_s} \|z - z_k\|^2 + \langle \widetilde{\nabla}_{k+1}, z \rangle + \psi(z) \big\}$;
    Option I: $y_{k+1} \leftarrow \arg\min_y \big\{ \frac{3L}{2} \|y - x_{k+1}\|^2 + \langle \widetilde{\nabla}_{k+1}, y \rangle + \psi(y) \big\}$;
    Option II: $y_{k+1} \leftarrow x_{k+1} + \tau_{1,s}(z_{k+1} - z_k)$ ◇ *we analyze only I but II also works*
  **end for**
  $\widetilde{x}^{s+1} \leftarrow \frac{1}{m} \sum_{j=1}^{m} y_{sm+j}$; ◇ *compute snapshot* $\widetilde{x}$
**end for**
**return** $\widetilde{x}^S$.

The only difference between $\texttt{Katyusha}^{\textsf{ns}}$ and `Katyusha` is that we choose $\tau_1 = \tau_{1,s} = \frac{2}{s+4}$ to be a parameter that depends on the epoch index $s$, and accordingly $\alpha = \alpha_s = \frac{1}{3L\tau_{1,s}}$. This should not be a big surprise because in accelerated full-gradient methods, the values $\tau_1$ and $\alpha$ also decrease (although with respect to $k$ rather than $s$) when there is no strong convexity [12]. We note that $\tau_1$ and $\tau_2$ remain constant throughout an epoch, and this could simplify the implementations.

We state the following convergence theorem for $\texttt{Katyusha}^{\textsf{ns}}$ and defer its proof to Appendix C.1. The proof also relies on the one-iteration inequality in Lemma 2.7, but requires telescoping such inequalities in a different manner as compared with Theorem 2.1.

**Theorem 4.1.** *If each $f_i(x)$ is convex, $L$-smooth in Problem (1.1) and $\psi(\cdot)$ is not necessarily strongly convex, then* $\texttt{Katyusha}^{\textsf{ns}}(x_0, S, L)$ *satisfies*

$$\mathbb{E}\big[F(\widetilde{x}^S)\big] - F(x^*) \leq O\Big(\frac{F(x_0) - F(x^*)}{S^2} + \frac{L\|x_0 - x^*\|^2}{mS^2}\Big)$$

*In other words, choosing $m = \Theta(n)$,* $\texttt{Katyusha}^{\textsf{ns}}$ *achieves an $\varepsilon$-additive error (i.e., $\mathbb{E}\big[F(\widetilde{x}^S)\big] - F(x^*) \leq \varepsilon$) using at most $O\Big(\frac{n\sqrt{F(x_0)-F(x^*)}}{\sqrt{\varepsilon}} + \frac{\sqrt{nL}\|x_0 - x^*\|}{\sqrt{\varepsilon}}\Big)$ iterations.*

*Remark* 4.2. $\texttt{Katyusha}^{\textsf{ns}}$ is a *direct, accelerated* solver for the non-SC case of Problem (1.1). It is illustrative to compare it with the convergence theorem of a *direct, non-accelerated* solver of the same setting. Below is the convergence theorem of SAGA after translating to our notations:

$$\text{SAGA:} \quad \mathbb{E}\big[F(x)\big] - F(x^*) \leq O\Big(\frac{F(x_0) - F(x^*)}{S} + \frac{L\|x_0 - x^*\|^2}{nS}\Big) \ .$$

It is clear from this comparison that $\texttt{Katyusha}^{\textsf{ns}}$ is a factor $S$ faster than non-accelerated methods such as SAGA, where $S = T/n$ if $T$ is the total number of stochastic iterations. This convergence can also be written in terms of the number of iterations which is $O\big(\frac{n(F(x_0)-F(x^*))}{\varepsilon} + \frac{L\|x_0-x^*\|^2}{\varepsilon}\big)$.

*Remark* 4.3. Theorem 4.1 appears worse than the reduction-based complexity in Corollary 3.7. This can be fixed by setting either the parameters $\tau_1$ or the epoch length $m$ in a more sophisticated way. Since it complicates the proofs and the notations we refrain from doing so in this version of the paper.[13] In practice, being a direct method, $\texttt{Katyusha}^{\mathsf{ns}}$ enjoys satisfactory performance.

# 5 Katyusha in the Mini-Batch Setting

We mentioned in earlier versions of this paper that our Katyusha method naturally generalizes to mini-batch (parallel) settings and non-uniform smoothness settings, but did not include a full proof. In this section, we carefully deal with both generalizations together.

**Mini-batch.** In each iteration $k$, instead of using a single $\nabla f_i(x_{k+1})$, one can

$$\text{use the average of } b \text{ stochastic gradients } \frac{1}{b}\sum_{i\in S_k}\nabla f_i(x_{k+1})$$

where $S_k$ is a random subset of $[n]$ with cardinality $b$. This average can be computed in a distributed manner using up to $b$ processors. This idea is known as *mini-batch* for stochastic gradient methods.

**Non-Uniform Smoothness.** Suppose in Problem (1.1),

$$\text{each } f_i(x) \text{ is } L_i\text{-smooth and } f(x) = \frac{1}{n}\sum_{i=1}^n f_i(x) \text{ is } L\text{-smooth.}$$

We denote by $\overline{L} = \frac{1}{n}\sum_{i=1}^n L_i$, and assume without loss of generality $L \le \overline{L} \le nL$. [14] We note that $\overline{L}$ can sometimes be indeed much greater than $L$.

*Remark* 5.1. In linear regression, we have $f_i(x) = \frac{1}{2}(\langle a_i, x\rangle - b_i)^2$. Denoting by $A = [a_1, \dots, a_n] \in \mathbb{R}^{d\times n}$, we have $L = \frac{1}{n}\lambda_{\max}(A^\top A)$ and $\overline{L} = \frac{1}{n}\|A\|_F^2$. If each entry of each $a_i$ is a random Gaussian $N(0,1)$, then $\overline{L}$ is around $d$ and $\overline{L}$ is around only $\Theta(1 + \frac{d}{n}) \ll \overline{L}$.

*Remark* 5.2. $L_i$ and $L$ only need to be upper bounds to the minimum smoothness parameters of $f_i(\cdot)$ and $f(\cdot)$ respectively. In practice, sometimes the minimum smoothness parameters for $f_i(x)$ is efficiently computable (such as for ERM problems).

## 5.1 Algorithmic Changes and Theorem Restatement

To simultaneously deal with mini-batch and non-uniform smoothness, we propose the following changes to `Katyusha`:

- *Change the epoch length from* $m = \Theta(n)$ *to* $m = \lceil \frac{n}{b} \rceil$.

  This is standard. In each iteration we need to compute $O(b)$ stochastic gradients; therefore every $\lceil \frac{n}{b} \rceil$ iterations, we can compute the full gradient once without hurting the total performance.

- *Define distribution* $\mathcal{D}$ *over* $[n]$ *to be choosing* $i \in [n]$ *with probability* $p_i \stackrel{\text{def}}{=} L_i / n\overline{L}$, *and define gradient estimator* $\widetilde{\nabla}_{k+1} \stackrel{\text{def}}{=} \nabla f(\widetilde{x}) + \frac{1}{b}\sum_{i\in S_k}\frac{1}{np_i}\big(\nabla f_i(x_{k+1}) - \nabla f_i(\widetilde{x})\big)$, *where* $S_k \subseteq [n]$ *is a multiset with* $b$ *elements each i.i.d. generated from* $\mathcal{D}$.

  This is standard, see for instance Prox-SVRG [57], and it is easy to verify $\mathbb{E}[\widetilde{\nabla}_{k+1}] = \nabla f(x_{k+1})$.

- *Change* $\tau_2$ *from* $\frac{1}{2}$ *to* $\min\big\{\frac{\overline{L}}{2Lb}, \frac{1}{2}\big\}$.

  Note that if $\overline{L} = L$ then we have $\tau_2 = \frac{1}{2b}$. In other words, the larger the mini-batch size, the smaller weight we want to give to Katyusha momentum. This should be intuitive.

[13] Recall that a similar issue has also happened in the non-accelerated world: the iteration complexity $O(\frac{n+L}{\varepsilon})$ in SAGA can be improved to $O(n\log\frac{1}{\varepsilon} + \frac{L}{\varepsilon})$ by doubling the epoch length across epochs [13]. Similar techniques can also be used to improve our result above.

[14] It is easy to verify (using triangle inequality) that $f(x) = \frac{1}{n}\sum_{i\in[n]} f_i(x)$ must be $\overline{L}$ smooth. Also, if $f(x)$ is $L$-smooth then each $f_i(x)$ must be $nL$ smooth (this can be checked via Hessian $\nabla^2 f_i(x) \preceq n\nabla^2 f(x)$ or similarly if $f$ is not twice-differentiable).

**Algorithm 3** `Katyusha1`$(x_0, S, \sigma, L, (L_1, \dots, L_n), b)$

1: $m \leftarrow \lceil n/b \rceil$ and $\overline{L} \leftarrow \frac{1}{n}(L_1 + \cdots + L_n)$; ◇ *$m$ is epoch length*
2: $\tau_2 \leftarrow \min\left\{\frac{\overline{L}}{2Lb}, \frac{1}{2}\right\}$; ◇ *if $\overline{L} = L$ then $\tau_2 = \frac{1}{2b}$ and $L_\diamond = L$*
3: **if** $L \le \frac{\overline{L}m}{b}$ **then**
4: $\quad\tau_1 \leftarrow \min\left\{\frac{\sqrt{8bm\sigma}}{\sqrt{3\overline{L}}}\tau_2, \tau_2\right\}$ and $L_\diamond \leftarrow \frac{\overline{L}}{2b\tau_2}$;
5: **else**
6: $\quad\tau_1 \leftarrow \min\left\{\frac{\sqrt{2\sigma}}{\sqrt{3L}}, \frac{1}{2m}\right\}$ and $L_\diamond \leftarrow L$;
7: **end if**
8: $\alpha \leftarrow \frac{1}{3\tau_1 L_\diamond}$; ◇ *parameters*
9: Let distribution $\mathcal{D}$ be to output $i \in [n]$ with probability $p_i \stackrel{\text{def}}{=} L_i/(n\overline{L})$.
10: $y_0 = z_0 = \widetilde{x}^0 \leftarrow x_0$; ◇ *initial vectors*
11: **for** $s \leftarrow 0$ **to** $S-1$ **do**
12: $\quad\mu^s \leftarrow \nabla f(\widetilde{x}^s)$; ◇ *compute the full gradient once every $m$ iterations*
13: $\quad$**for** $j \leftarrow 0$ **to** $m-1$ **do**
14: $\quad\quad k \leftarrow (sm) + j$;
15: $\quad\quad x_{k+1} \leftarrow \tau_1 z_k + \tau_2 \widetilde{x}^s + (1 - \tau_1 - \tau_2) y_k$;
16: $\quad\quad S_k \leftarrow b$ independent copies of $i$ from $\mathcal{D}$ with replacement.
17: $\quad\quad \widetilde{\nabla}_{k+1} \leftarrow \mu^s + \frac{1}{b}\sum_{i \in S_k} \frac{1}{np_i}\big(\nabla f_i(x_{k+1}) - \nabla f_i(\widetilde{x}^s)\big)$;
18: $\quad\quad z_{k+1} = \arg\min_z \left\{\frac{1}{2\alpha}\|z - z_k\|^2 + \langle \widetilde{\nabla}_{k+1}, z\rangle + \psi(z)\right\}$;
19: $\quad\quad$Option I: $y_{k+1} \leftarrow \arg\min_y \left\{\frac{3L_\diamond}{2}\|y - x_{k+1}\|^2 + \langle \widetilde{\nabla}_{k+1}, y\rangle + \psi(y)\right\}$;
20: $\quad\quad$Option II: $y_{k+1} \leftarrow x_{k+1} + \tau_1(z_{k+1} - z_k)$ ◇ *we analyze only I but II also works*
21: $\quad$**end for**
22: $\quad\widetilde{x}^{s+1} \leftarrow \big(\sum_{j=0}^{m-1} \theta^j\big)^{-1} \cdot \big(\sum_{j=0}^{m-1} \theta^j \cdot y_{sm+j+1}\big)$; ◇ *where $\theta = 1 + \min\{\alpha\sigma, \frac{1}{4m}\}$*
23: **end for**
24: **return** $x^{\mathsf{out}} \leftarrow \frac{\tau_2 m \widetilde{x}^S + (1-\tau_1-\tau_2) y_{Sm}}{\tau_2 m + (1-\tau_1-\tau_2)}$.

- *Change $L$ in gradient descent step (Line 19) to some other $L_\diamond \ge L$, and define $\alpha = \frac{1}{3\tau_1 L_\diamond}$ instead.*

  In most cases (e.g., when $\overline{L} = L$ or $L \ge \overline{L}m/b$) we choose $L_\diamond = L$. Otherwise, we let $L_\diamond = \frac{\overline{L}}{2b\tau_2} \ge L$.

- *Change $\tau_1$ to be $\tau_1 = \min\left\{\frac{\sqrt{8bm\sigma}}{\sqrt{3\overline{L}}}\tau_2, \tau_2\right\}$ if $L \le \overline{L}m/b$ or $\tau_1 = \min\left\{\frac{\sqrt{2\sigma}}{\sqrt{3L}}, \frac{1}{2m}\right\}$ if $L > \overline{L}m/b$.*

  This corresponds to a phase-transition behavior of `Katyusha1` (see Remark 5.6 later). Intuitively, when $L \le \overline{L}m/b$ then we are in a mini-batch phase; when $L > \overline{L}m/b$ we are in a full-batch phase.

- *Due to technical reasons, we define $\widetilde{x}^s$ as a slightly different weighted average (Line 22) and output $x^{\mathsf{out}}$ which is a weighted combination of $\widetilde{x}^S$ and $y_{Sm}$ as opposed to simply $\widetilde{x}^S$ (Line 24).*

We emphasize here that some of these changes are not necessary for instance in the special case of $\overline{L} = L$, but to state the strongest theorem, we have to include all such changes. It is a simple exercise to verify that, if $\overline{L} = L$ and $b = 1$, then up to only constant factors in the parameters, `Katyusha1` is exactly identical to `Katyusha`. We have the following main theorem for `Katyusha1`:

*Remark* 5.3. If $\overline{L}$ is equal to $L$, then one can simply set $\tau_2 = \frac{1}{2b}$ and $L_\diamond = L$ in `Katyusha1`.

**Theorem 5.4.** *If each $f_i(x)$ is convex and $L_i$-smooth, $f(x)$ is $L$-smooth, $\psi(x)$ is $\sigma$-strongly convex in Problem (1.1), then for any $b \in [n]$, $x^{\mathsf{out}} = \texttt{Katyusha1}(x_0, S, \sigma, L, (L_1, \dots, L_n), b)$ satisfies*

$$\mathbb{E}\big[F(x^{\mathsf{out}})\big] - F(x^*) \leq \begin{cases} O\Big(\Big(1+\sqrt{b\sigma/(6\overline{L}m)}\Big)^{-Sm}\Big) \cdot \big(F(x_0) - F(x^*)\big), & \text{if } \frac{m\sigma b}{\overline{L}} \leq \frac{3}{8} \text{ and } L \leq \frac{\overline{L}m}{b}; \\ O\Big(\big(1+\sqrt{\sigma/(6L)}\big)^{-Sm}\Big) \cdot \big(F(x_0) - F(x^*)\big), & \text{if } \frac{m^2\sigma}{L} \leq \frac{3}{8} \text{ and } L > \frac{\overline{L}m}{b}; \\ O\big(1.25^{-S}\big) \cdot \big(F(x_0) - F(x^*)\big), & \text{otherwise.} \end{cases}$$

*In other words, choosing $m = \lceil n/b \rceil$, `Katyusha` achieves an $\varepsilon$-additive error (that is, $\mathbb{E}\big[F(x^{\mathsf{out}})\big] - F(x^*) \leq \varepsilon$) using at most*

$$S \cdot n = O\Big(\big(n + b\sqrt{L/\sigma} + \sqrt{n\overline{L}/\sigma}\big) \cdot \log \frac{F(x_0) - F(x^*)}{\varepsilon}\Big)$$

*stochastic gradient computations.*

### 5.2 Observations and Remarks

We explain the significance of Theorem 5.4 below. We use *total work* to refer to the total number of stochastic gradient computations, and *iteration complexity* (also known as parallel depth) to refer to the total number of iterations.

**Parallel Linear Speedup.** The total work of `Katyusha1` stays the same when $b \leq (n\overline{L}/L)^{1/2} \in \big[\sqrt{n}, n\big]$. This means, at least for all $b \in \{1, 2, \dots, \lceil\sqrt{n}\rceil\}$, our `Katyusha1` has the same total work.

`Katyusha1` can be distributed to $b \leq \sqrt{n}$ machines with a parallel speed-up factor $b$ (known as linear speedup if ignoring communication overhead.)

In contrast, for instance in the special case of $\overline{L} = L$, if no additional assumption is made:

- Mini-batch SVRG requires $\widetilde{O}\big(n + \frac{bL}{\sigma}\big)$ total work.

  Therefore, if SVRG is distributed to $b$ machines, the total work is increased by a factor of $b$, and the parallel speed-up factor is 1 (i.e., no speed up).

- Catalyst on top of mini-batch SVRG requires $\widetilde{O}\big(n + \frac{\sqrt{bLn}}{\sqrt{\sigma}}\big)$ total work.

  Therefore, if Catalyst is distributed to $b$ machines, the total work is increased by a factor $\sqrt{b}$, and the parallel speed-up factor is $\sqrt{b}$ only.

*Remark* 5.5. When preparing the journal revision (i.e., version 5), we found out at least in the case $\overline{L} = L$, Murata and Suzuki [41] independently obtained similar "parallel linear speedup." Their analysis is also primal-only. In contrast, it is unclear if any dual-based method (such as RPDG [33], AccSDCA [54], APCG [37], or SPDC [60]) enjoys such linear-speedup. We are aware of an attempt to tackle SPDC in the parallel setting [55], but it does not lead to a parallel linear speedup.

**Outperforming Nesterov's Method.** At first glance, if $b = n$, `Katyusha1` runs in $\widetilde{O}\big((L/\sigma)^{1/2}\big)$ iterations. This matches the total work and iteration complexity of Nesterov's accelerated gradient method [12, 42, 43], and does not depend on the possibly larger parameter $\overline{L}$.

More interestingly, to achieve the *same* iteration complexity $\widetilde{O}\big((L/\sigma)^{1/2}\big)$ as Nesterov's method, `Katyusha1` only needs to compute $b = (n\overline{L}/L)^{1/2}$ stochastic gradients $\nabla f_i(\cdot)$ per iteration in average. This can be much faster than computing $\nabla f(\cdot)$, thus giving a much smaller total work.

**Phase Transition between Mini-Batch and Full-Batch.** Theorem 5.4 indicates a phase transition of `Katyusha1` at the point $b_0 = (n\overline{L}/L)^{1/2}$.

- If $b \leq b_0$, we say `Katyusha1` is in the *mini-batch phase* and the total work is $\widetilde{O}\big(n + \sqrt{n\overline{L}/\sigma}\big)$, independent of $b$.
- If $b > b_0$, we say `Katyusha1` is in the *full-batch phase*, and the total work is $\widetilde{O}\big(n + b\sqrt{L/\sigma}\big)$, so essentially linearly-scales with $b$ and matches that of Nesterov's method when $b = n$.

*Remark* 5.6. We set different values for $\tau_1$ and $L_\diamond$ in the mini-batch phase and full-batch phase respectively (see Line 3). From the final complexities above, it should not be surprising that $\tau_1$ depends on $\overline{L}$ but not $L$ in the mini-batch phase, and depends on $L$ but not $\overline{L}$ in the full-batch phase. In addition, one can even tune the parameters so that it suffices for `Katyusha` to output $\widetilde{x}^S$ in the mini-batch phase and $y_{Sm}$ in the full-batch phase; we did not do so and simply choose to output $x^{\mathsf{out}}$ which is a convex combination of $\widetilde{x}^S$ and $y_{Sm}$.

*Remark* 5.7. In the simple case $\overline{L} = L$, Nitanda [46] obtained a total work $\widetilde{O}\big(n + \frac{n-b}{n-1}\frac{L}{\sigma} + b\sqrt{L/\sigma}\big)$ (see their page 7), which also implies a phase transition for $b$. However, this result is no better than ours for *all* parameters $b \in \{1, 2, \dots, n\}$. In fact, in terms of total work, Nitanda's result is even slower than SVRG when $b \leq n/2$, and slower than Nesterov's method when $b > n/2$.

### 5.3 Corollaries on Non-Smooth or Non-SC Problems

In the same way as Section 3, we can apply the reductions from [6] to convert the performance of Theorem 5.4 to non-smooth or non-strongly convex settings. We state the corollaries below:

**Corollary 5.8.** *If each $f_i(x)$ is convex and $L_i$-smooth, $f(x)$ is $L$-smooth, $\psi(\cdot)$ is not necessarily strongly convex in Problem (1.1), then for any $b \in [n]$, by applying* AdaptReg *on* `Katyusha1` *with a starting vector $x_0$, we obtain an output $x$ satisfying $\mathbb{E}[F(x)] - F(x^*) \leq \varepsilon$ in at most*

$$O\Big(n \log \tfrac{F(x_0)-F(x^*)}{\varepsilon} + \tfrac{b\sqrt{L}\cdot\|x_0-x^*\|}{\sqrt{\varepsilon}} + \tfrac{\sqrt{n\overline{L}}\cdot\|x_0-x^*\|}{\sqrt{\varepsilon}}\Big) \text{ stochastic gradient computations.}$$

**Corollary 5.9.** *If each $f_i(x)$ is $\sqrt{G_i}$-Lipschitz continuous and $\psi(x)$ is $\sigma$-SC in Problem (1.3), then for any $b \in [n]$, by applying* AdaptSmooth *on* `Katyusha1` *with a starting vector $x_0$, we obtain an output $x$ satisfying $\mathbb{E}[F(x)] - F(x^*) \leq \varepsilon$ in at most*

$$O\Big(n \log \tfrac{F(x_0)-F(x^*)}{\varepsilon} + \tfrac{b\sqrt{G}}{\sqrt{\sigma\varepsilon}} + \tfrac{\sqrt{n\overline{G}}}{\sqrt{\sigma\varepsilon}}\Big) \text{ stochastic gradient computations.}$$

**Corollary 5.10.** *If each $f_i(x)$ is $\sqrt{G_i}$-Lipschitz continuous and $\psi(x)$ is not necessarily strongly convex in Problem (1.3), then for any $b \in [n]$, by applying* JointAdaptRegSmooth *on* `Katyusha1` *with a starting vector $x_0$, we obtain an output $x$ satisfying $\mathbb{E}[F(x)] - F(x^*) \leq \varepsilon$ in at most*

$$O\Big(n \log \tfrac{F(x_0)-F(x^*)}{\varepsilon} + \tfrac{b\,\overline{G}\|x_0-x^*\|}{\varepsilon} + \tfrac{\sqrt{n}\,\overline{G}\|x_0-x^*\|}{\varepsilon}\Big) \text{ stochastic gradient computations.}$$

## 6 Katyusha in the Non-Euclidean Norm Setting

In this section, we show that `Katyusha` and `Katyusha`$^{\mathsf{ns}}$ naturally extend to settings where the smoothness definition is with respect to a non-Euclidean norm.

**Non-Euclidean Norm Smoothness.** We consider smoothness (and strongly convexity) with respect to an arbitrary norm $\|\cdot\|$ in domain $Q \stackrel{\text{def}}{=} \{x \in \mathbb{R}^d : \psi(x) < +\infty\}$. Symbolically, we say

- $f$ is $\sigma$-strongly convex w.r.t. $\|\cdot\|$ if $\forall x, y \in Q$, it satisfies $f(y) \geq f(x) + \langle \nabla f(x), y-x\rangle + \frac{\sigma}{2}\|x-y\|^2$;

- $f$ is $L$-smooth w.r.t. $\|\cdot\|$ if $\forall x,y\in Q$, it satisfies $\|\nabla f(x)-\nabla f(y)\|_* \le L\|x-y\|$.[15]

Above, $\|\cdot\|_* \stackrel{\text{def}}{=} \max\{\langle \xi, x\rangle : \|x\|\le 1\}$ is the dual norm of $\|\cdot\|$. For instance, $\ell_p$ norm is dual to $\ell_q$ norm if $\frac{1}{p}+\frac{1}{q}=1$. Some famous problems have better smoothness parameters when non-Euclidean norms are adopted, see the discussions in [12].

**Bregman Divergence.** Following the traditions in the non-Euclidean norm setting [12], we

- select a *distance generating function* $w(\cdot)$ that is 1-strongly convex w.r.t. $\|\cdot\|$, and[16]
- define the *Bregman divergence function* $V_x(y) \stackrel{\text{def}}{=} w(y)-w(x)-\langle \nabla w(x), y-x\rangle$.

The final algorithms and proofs will be described using $V_x(y)$ and $w(x)$.

**Generalized Strong Convexity of $\psi(\cdot)$.** We require $\psi(\cdot)$ to be $\sigma$-strongly convexity with respect to function $V_x(y)$ rather than the $\|\cdot\|$ norm; or symbolically,

$$\psi(y) \ge \psi(x) + \langle \nabla\psi(x), y-x\rangle + \sigma V_x(y) \enspace .$$

(For instance, this is satisfied if $\omega(y) \stackrel{\text{def}}{=} \frac{1}{\sigma}\psi(y)$.) This is known as the "generalized strong convexity" [50] and is necessary for any linear-convergence result in the SC setting. Of course, in the non-SC setting, we do not require any (general or not) strong convexity for $\psi(\cdot)$.

### 6.1 Algorithm Changes and Theorem Restatements

Suppose each $f_i(x)$ is $L_i$-smooth with respect to norm $\|\cdot\|$, and a Bregman divergence function $V_x(y)$ is given. We perform the following changes to the algorithms:

- In Line 9 of `Katyusha` (resp. Line 10 of `Katyusha`$^{\mathsf{ns}}$), we choose $i$ with probability proportional to $L_i$ instead of uniformly at random.
- In Line 10 of `Katyusha` (resp. Line 11 of `Katyusha`$^{\mathsf{ns}}$), we change the arg min to be its non-Euclidean norm variant [12]: $z_{k+1} = \arg\min_z \left\{\frac{1}{\alpha}V_{z_k}(z) + \langle \widetilde{\nabla}_{k+1}, z\rangle + \psi(z)\right\}$
- We forbidden Option II and use Option I only (but *without* replacing $\|y-x_{k+1}\|^2$ with $V_{x_{k+1}}(y)$).

Interested readers can find discussions regarding why such changes are natural in [12]. We call the resulting algorithms `Katyusha2` and `Katyusha2`$^{\mathsf{ns}}$, and include them in Appendix E for completeness' sake. We state our final theorems below (recall $\overline{L} = \frac{1}{n}\sum_{i=1}^n L_i$).

**Theorem 6.1** (ext. of Theorem 2.1)**.** *If each $f_i(x)$ is convex and $L_i$-smooth with respect to some norm $\|\cdot\|$, $V_x(y)$ is a Bregman divergence function for $\|\cdot\|$, and $\psi(x)$ is $\sigma$-strongly convex with respect to $V_x(y)$, then* `Katyusha2`$(x_0, S, \sigma, (L_1,\dots,L_n))$ *satisfies*

$$\mathbb{E}\big[F(\widetilde{x}^S)\big] - F(x^*) \le \begin{cases} O\Big(\Big(1+\sqrt{\sigma/(9\overline{L}m)}\Big)^{-Sm}\Big)\cdot\big(F(x_0)-F(x^*)\big), & \text{if } m\sigma/\overline{L} \le \frac{9}{4}; \\ O\big(1.5^{-S}\big)\cdot\big(F(x_0)-F(x^*)\big), & \text{if } m\sigma/\overline{L} > \frac{9}{4}. \end{cases}$$

*In other words, choosing $m=\Theta(n)$,* `Katyusha2` *achieves an $\varepsilon$-additive error (i.e., $\mathbb{E}\big[F(\widetilde{x}^S)\big] - F(x^*) \le \varepsilon$) using at most $O\big(\big(n+\sqrt{n\overline{L}/\sigma}\big)\cdot\log\frac{F(x_0)-F(x^*)}{\varepsilon}\big)$ iterations.*

[15] This definition has another equivalent form: $\forall x,y\in Q$, it satisfies $f(y)\le f(x)+\langle\nabla f(x),y-x\rangle+\frac{L}{2}\|y-x\|^2$.

[16] For instance, if $Q=\mathbb{R}^d$ and $\|\cdot\|_p$ is the $\ell_p$ norm for some $p\in(1,2]$, one can choose $w(x)=\frac{1}{2(p-1)}\|x\|_p^2$; if $Q=\{x\in\mathbb{R}^d : \sum_i x_i = 1\}$ is the probability space and $\|\cdot\|_1$ is the $\ell_1$ norm, one can choose $w(x)=\sum_i x_i\log x_i$.

**Theorem 6.2** (ext. of Theorem 4.1)**.** *If each $f_i(x)$ is convex and $L_i$-smooth with respect to some norm $\|\cdot\|$, $V_x(y)$ is a Bregman divergence function for $\|\cdot\|$, and $\psi(\cdot)$ is not necessarily strongly convex, then* $\texttt{Katyusha2}^{\textsf{ns}}(x_0, S, (L_1, \dots, L_n))$ *satisfies*

$$\mathbb{E}\big[F(\widetilde{x}^S)\big] - F(x^*) \le O\Big(\frac{F(x_0) - F(x^*)}{S^2} + \frac{\overline{L}V_{x_0}(x^*)}{nS^2}\Big) \ .$$

*In other words,* $\texttt{Katyusha2}^{\textsf{ns}}$ *achieves an $\varepsilon$-additive error (i.e., $\mathbb{E}\big[F(\widetilde{x}^S)\big] - F(x^*) \le \varepsilon$) using at most $O\Big(\frac{n\sqrt{F(x_0)-F(x^*)}}{\sqrt{\varepsilon}} + \frac{\sqrt{nLV_{x_0}(x^*)}}{\sqrt{\varepsilon}}\Big)$ iterations.*

The proofs of Theorem 6.1 and Theorem 6.2 follow exactly the same proof structures of Theorem 2.1 and Theorem 4.1, so we include them only in Appendix E.

### 6.2 Remarks

We highlight one main difference between the proof of `Katyusha2` and that of `Katyusha`: if $\xi$ is a random vector and $\|\cdot\|$ is an arbitrary norm, we do not necessarily have $\mathbb{E}[\|\xi - \mathbb{E}[\xi]\|_*^2] \le \mathbb{E}[\|\xi\|_*^2]$. Therefore, we only used $\mathbb{E}[\|\xi - \mathbb{E}[\xi]\|_*^2] \le 2\mathbb{E}[\|\xi\|_*^2] + 2\|\mathbb{E}[\xi]\|_*^2$ (see Lemma E.2) and this loses a constant factor in some parameters. (For instance, $\alpha$ now becomes $\frac{1}{9\tau_1\overline{L}}$ as opposed to $\frac{1}{3\tau_1 L}$).

More interestingly, one may ask how our revised algorithms `Katyusha2` or $\texttt{Katyusha2}^{\textsf{ns}}$ perform in the mini-batch setting (just like we have studied in Section 5 for the Euclidean case). We are optimistic here, but unfortunately do not have a clean worst-case statement for how much speed-up we can get. The underlying reason is that, if $\mathcal{D}$ is a distribution for vectors, $\mu = \mathbb{E}_{\xi\sim\mathcal{D}}[\xi]$ is its expectation, and $\xi_1, \dots, \xi_b$ are $b$ i.i.d. samples from $\mathcal{D}$, then letting $\overline{\xi} = \frac{1}{b}(\xi_1 + \cdots + \xi_b)$, we do not necessarily have $\mathbb{E}[\|\overline{\xi} - \mu\|_*^2] \le \frac{1}{b}\mathbb{E}_{\xi\sim\mathcal{D}}[\|\xi - \mu\|_*^2]$. In other words, using a mini-batch version of the gradient estimator, the "variance" with respect to an arbitrary norm may not necessarily go down by a factor of $b$. For such reason, in the mini-batch setting, the best total work we can cleanly state, say for `Katyusha2` in the SC setting, is only $O\Big(\big(n + \sqrt{bn\overline{L}/\sigma}\big) \cdot \log\frac{F(x_0)-F(x^*)}{\varepsilon}\Big)$ .

## 7 Empirical Evaluations

We conclude this paper with empirical evaluations to our theoretical speed-ups. We work on Lasso and ridge regressions (with regularizer $\frac{\lambda}{2}\|x\|^2$ for ridge and regularizer $\lambda\|x\|_1$ for Lasso) on the following six datasets: adult, web, mnist, rcv1, covtype, sensit. We defer dataset and implementation details to Appendix B.

**Algorithms and Parameter Tuning.** We have implemented the following algorithms, all with mini-batch size 1 for this version of the paper:

- SVRG [29] with default epoch length $m = 2n$. We tune only *one parameter*: the learning rate.
- `Katyusha` for ridge and $\texttt{Katyusha}^{\textsf{ns}}$ for Lasso. We tune only *one parameter*: the learning rate.
- SAGA [22]. We tune only *one parameter*: the learning rate.
- Catalyst [36] on top of SVRG. We tune *three parameters*: SVRG's learning rate, Catalyst's learning rate, as well as the regularizer weight in the Catalyst reduction.
- APCG [37]. We tune the learning rate. For Lasso, we also tune the $\ell_2$ regularizer weight.
- APCG+AdaptReg (Lasso only). Since APCG intrinsically require an $\ell_2$ regularizer to be added on Lasso, we apply AdaptReg from [6] to adaptively learn this regularizer and improve APCG's performance. Two parameters to be tuned: APCG's learning rate and $\sigma_0$ in AdaptReg.

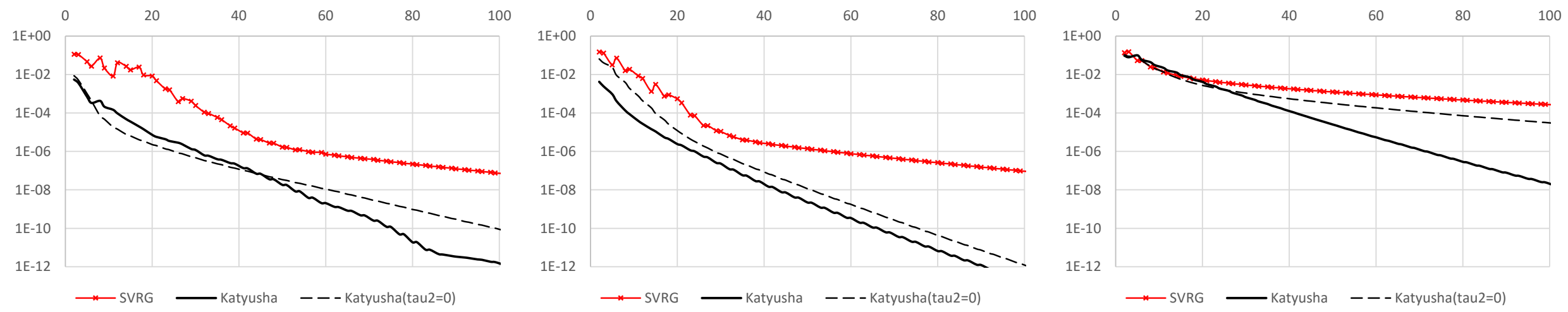

(a) web, ridge $\lambda = 10^{-6}$ (b) mnist, ridge $\lambda = 10^{-6}$ (c) rcv1, ridge $\lambda = 10^{-6}$

Figure 1: Comparing SVRG vs. `Katyusha` vs. `Katyusha` with $\tau_2 = 0$.

All of the parameters were equally, fairly, and automatically tuned by our code base. For interested readers, we discuss more details in Appendix B.

We emphasize that `Katyusha` is *as simple as SAGA or SVRG in terms of parameter tuning.* In contrast, APCG for Lasso requires two parameters to be tuned, and Catalyst requires three. [35]

**Performance Plots.** Following the tradition of ERM experiments, we use the number of "passes" of the dataset as the $x$-axis. Letting $n$ be the number of feature vectors, each new stochastic gradient computation $\nabla f_i(\cdot)$ counts as $1/n$ pass, and a full gradient computation $\nabla f(\cdot)$ counts as 1 pass.

The $y$-axis in our plots represents the training objective distance to the minimum. Since we aim to evaluate our theoretical finding, we did not include the test error. We emphasize that it is practically also crucial to study high-accuracy regimes (such as objective distance $\leq 10^{-7}$). This is because nowadays there is an increasing number of methods that reduce large-scale machine learning tasks to multiple black-box calls to ERM solvers [8, 9, 25]. In all such applications, due to error blowups between oracle calls, the ERM solver is required to be *very accurate in training error.*

## 7.1 Effectiveness of Katyusha Momentum

In our `Katyusha` method, $\tau_1$ controls to the classical Nesterov's momentum and $\tau_2$ controls our newly introduced Katyusha momentum. We find in our theory that setting $\tau_2 = 1/2$ is a good choice so we universally set it to be $1/2$ without tuning in all our experiments. (Of course, if time permits, tuning $\tau_2$ could only help in performance.)

Before this paper, researchers have tried heuristics that is to add Nesterov's momentums directly to stochastic gradient methods [46], and this corresponds to setting $\tau_2 = 0$ in `Katyusha`. In Figure 1, we compare `Katyusha` with $\tau_2 = 1/2$ and $\tau_2 = 0$ in order to illustrate the importance and effectiveness of our Katyusha momentum.

We conclude that the old heuristics (i.e., $\tau_2 = 0$) sometimes indeed make the method faster after careful parameter tuning. However, for certain tasks such as Figure 1(c), without Katyusha momentum the algorithm does not even enjoy an accelerated convergence rate.

## 7.2 Performance Comparison Across Algorithms

For each of the six datasets and each objective (ridge or lasso), we experiment on three different magnitudes of regularizer weights.[17] This totals 36 performance charts, and we include them in

[17] We choose three values $\lambda$ that are powers of 10 and around $10/n, 1/n, 1/10n$. This range can be verified to contain the best regularization weights using cross validation.

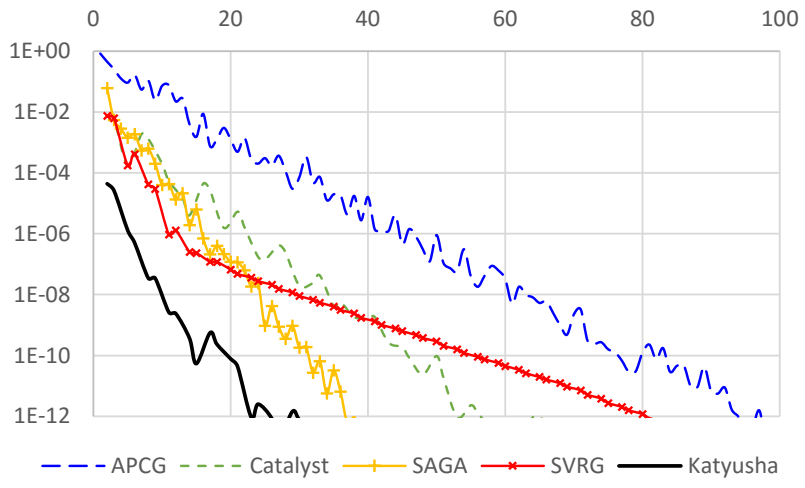

(a) covtype, ridge $\lambda = 10^{-6}$ (small)

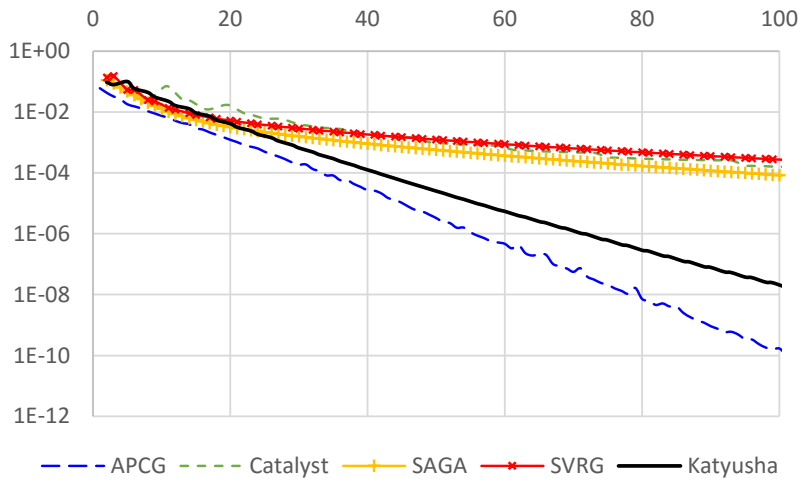

(b) rcv1, ridge $\lambda = 10^{-7}$ (small)

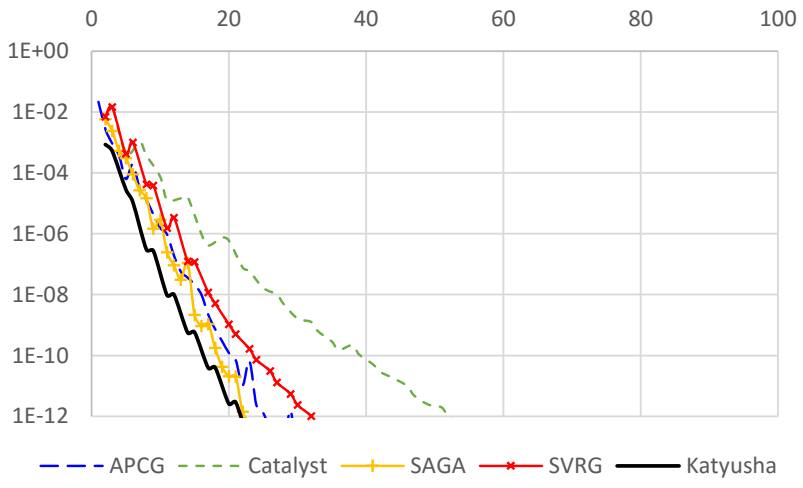

(c) adult, ridge $\lambda = 10^{-4}$ (large)

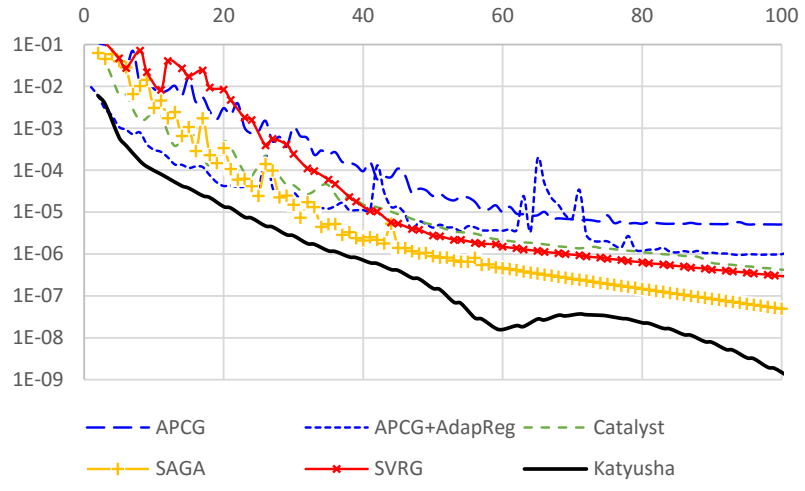

(d) web, ridge $\lambda = 10^{-6}$ (small)

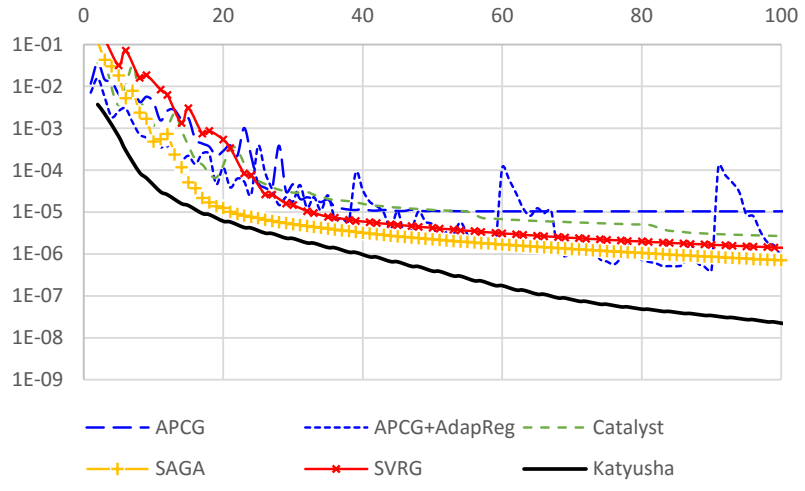

(e) mnist, ridge $\lambda = 10^{-6}$ (small)

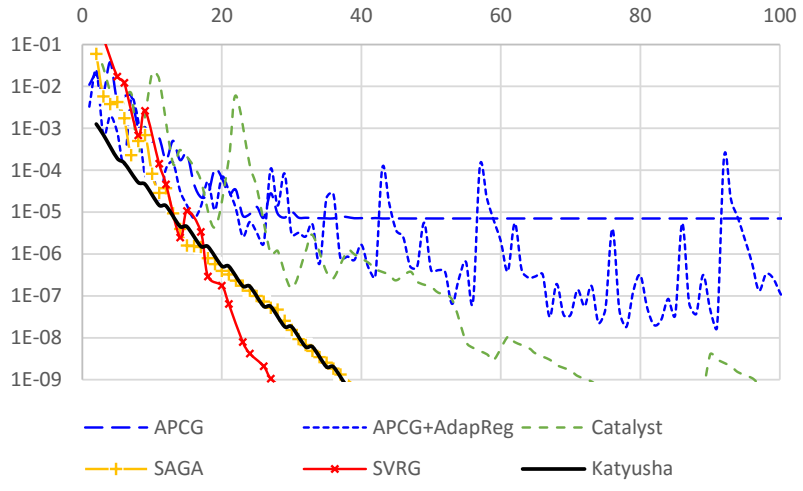

(f) sensit, ridge $\lambda = 10^{-4}$ (large)

Figure 2: Some representative performance charts where $\lambda$ is the regularizer weight. See Figure 3 and Figure 4 in the appendix for the full plots.

full at the end of this paper. For the sake of cleanness, in Figure 2 we select 6 representative charts for ridge regression and make the following observations.

- Accelerated methods are more powerful when the regularizer weights are small (cf. [16, 37, 54]). For instance, Figure 2(c) and 2(f) are for large values of $\lambda$ and `Katyusha` performs relatively the same as compared with SVRG / SAGA; however, `Katyusha` significantly outperforms SVRG / SAGA for small values of $\lambda$, see for instance Figure 2(b) and 2(e).
- `Katyusha` almost always either outperform or equal-perform its competitors. The only notable place it gets outperformed is by SVRG (see Figure 2(f)); however, this performance gap cannot be large because `Katyusha` is capable of recovering SVRG if $\tau_1 = \tau_2 = 0$.[18]
- Catalyst does not work as beautiful as its theory in high-accuracy regimes, even though we have carefully tuned parameters $\alpha_0$ and $\kappa$ in Catalyst in addition to its learning rate. Indeed, in Figure 2(a), 2(c) and 2(f) Catalyst (which is a reduction on SVRG) is outperformed by SVRG.
- APCG performs poorly on all Lasso tasks (cf. Figure 2(d), 2(e), 2(f)) because it is not designed for non-SC objectives. The reduction in [6] helps to fix this issue, but not by a lot.
- APCG can sometimes be largely dominated by SVRG or SAGA (cf. Figure 2(f)): this is because for datasets such as sensit, dual-based methods (such as APCG) cannot make use of the implicity local strong convexity in the objective. In such cases, `Katyusha` is not lost to SVRG or SAGA.

[18] The only reason `Katyusha` does not match the performance of SVRG in Figure 2(f) is because we have not tuned parameter $\tau_2$. If we also tune $\tau_2$ for the best performance, `Katyusha` shall no longer be outperformed by SVRG. In any case, it is not really necessary to tune $\tau_2$ because the performance of `Katyusha` is already superb.

# 8 Conclusion

The Katyusha momentum technique introduced in this paper gives rise to accelerated convergence rates even in the stochastic setting. For many classes of the problems, such convergence rates are the first to match the theoretical lower bounds [56]. The algorithms generated by Katyusha momentum are simple yet highly practical and parallelizable.

More importantly, this new technique has the potential to enrich our understanding of accelerated methods in a broader sense. Currently, although acceleration methods are becoming more and more important to the field of computer science, they are still often regarded as "analytical tricks" [19, 30] and lacking complete theoretical understanding. The Katyusha momentum presented in this paper, however, adds a new level of decoration on top of the classical Nesterov momentum. This decoration is shown valuable for stochastic problems in this paper, but may also lead to future applications as well. In general, the author hopes that the technique and analysis in this paper could facilitate more studies in this field and thus become a stepping stone towards the ultimate goal of unveiling the mystery of acceleration.

# Appendix

## A Other Related Works

For smooth convex minimization problems, (full) gradient descent converges at a rate $\frac{L}{\varepsilon}$—or $\frac{L}{\sigma}\log\frac{1}{\varepsilon}$ if the objective is $\sigma$-strongly convex. This is not optimal among the class of first-order methods. Nesterov showed that the optimal rate should be $\frac{\sqrt{L}}{\sqrt{\varepsilon}}$—or $\frac{\sqrt{L}}{\sqrt{\sigma}}\log\frac{1}{\varepsilon}$ if the objective is $\sigma$-strongly convex— and this was achieved by his celebrated accelerated (full) gradient descent method [42].

**Sum-of-Nonconvex Optimization.** One important generalization of Problem (1.1) is the case when the functions $f_i(x)$ are non-convex but their average $f(x) = \frac{1}{n}\sum_{i=1}^n f_i(x)$ is convex. Solvers for this "sum-of-nonconvex" setting can be applied to PCA/SVD, online eigenvector, general non-convex optimization, and more. See [3] and the references therein.

Variance reduction was first introduced to solve this problem by Shalev-Shwartz [51], and APPA/Catalyst also accelerates SVRG for this problem [26]. One can similarly ask whether one can design a *directly accelerated* method for this more general problem, and this was achieved by the `KatyushaX` method in [3]. It is a sister paper to us but uses very different sets of techniques.

**Online Stochastic Optimization.** Some literatures also focus on the (more general) online variant of Problem (1.1), that is for $n$ being sufficiently large so that the algorithm cannot access $f(x)$ or compute its full gradient $\nabla f(x)$. In this regime, without additional assumption, the optimal convergence rate is $1/\varepsilon^2$ (or $1/\varepsilon$ if the function is strongly convex). This is obtained by SGD and its hybrid variants [28, 32].

**Coordinate Descent.** Another way to define gradient estimator is to set $\widetilde{\nabla}_k = d\nabla_j f(x_k)$ where $j$ is a random coordinate. This is *(randomized) coordinate descent* as opposed to stochastic gradient descent. Designing accelerated methods for coordinate descent is significantly easier than designing that for stochastic gradient descent, and has indeed been done in many previous results including [16, 37, 38, 45].[19] The fastest rate is $O\big(\sum_i \sqrt{L_i/\varepsilon}\big)$ where parameters $L_i$ correspond to the coordi-nate smoothness of $f(x)$ [16]. Coordinate descent *cannot* be applied to solve Problem (1.1) because in our *stochastic* setting, only one copy $\nabla f_i(\cdot)$ is computed in every iteration.

[19]The reason behind it can be understood as follows. If a function $f(\cdot)$ is $L$ smooth with respect to coordinate $j$, then a coordinate descent step $x' \leftarrow x - \frac{1}{L}\nabla_j f(x)\mathbf{e}_j$ always decreases the objective, i.e., $f(x + \frac{1}{L}\nabla_j f(x)\mathbf{e}_j) < f(x)$. In contrast, this is *false* for stochastic gradient descent, because $f(x_k - \eta\widetilde{\nabla}_k)$ may be even larger than $f(x_k)$.

**Hybrid Stochastic Methods.** Several recent results study hybrid methods with convergence rates that are generally *non-accelerated* and only accelerated in *extreme cases.*

- Hu et al. [28], Lan [32] obtained iteration complexity of the form $O(L/\sqrt{\varepsilon}+\sigma/\varepsilon^2)$ in the presence of stochastic gradient with variance $\sigma$. These results can be interpreted as follows, if $\sigma$ is very small, then one can directly apply Nesterov's accelerated gradient method and achieve $O(L/\sqrt{\varepsilon})$; or if $\sigma$ is large then they match the SGD iteration complexity $O(\sigma/\varepsilon^2)$. For Problem (1.1), these algorithms do not give faster running time than `Katyusha` unless $\sigma$ is very small.[20]
- Nitanda [46] adds momentum to the non-accelerated variance-reduction method in a naive manner. It corresponds to this paper but *without* Katyusha momentum (i.e., $\tau_2 = 0$). The theoretical running time of [46] is always slower than this paper, and even slower than the best between SVRG and Nesterov's method. See Remark 5.7.[21] We have included an experiment in Section 7.1 to illustrate why Katyusha momentum is necessary.

**Linear Coupling.** Allen-Zhu and Orecchia proposed a framework called *linear coupling* that facilitates the design of accelerated gradient methods [12]. The simplest use of linear coupling can reconstruct Nesterov's accelerated full-gradient method [12], or to provide faster coordinate descent [16]. More careful use of linear coupling can also give accelerated methods for non-smooth problems (such as positive LP [10, 11], positive SDP [14], matrix scaling [17]) or for general non-convex problems [5]. This present paper falls into this linear-coupling framework, but our Katyusha momentum technique was not present in any of these cited results.

**Acceleration in Nonconvex Optimization.** One can also ask how does acceleration help in non-convex optimization? This is a new area with active research going on.

In the deterministic setting, under standard Lipschitz smoothness, gradient descent finds a point $x$ with $\|\nabla f(x)\| \le \varepsilon$ in $O(\varepsilon^{-2})$ iterations [43], and acceleration is not known to help. If second-order Lipschitz smoothness is added, then one can use momentum to non-trivially improve the rate to $O(\varepsilon^{-1.75})$ [1, 20].

In the finite-sum stochastic setting, gradient descent finds a point $x$ with $\|\nabla f(x)\| \le \varepsilon$ in $T = O(n\varepsilon^{-2})$ stochastic gradient computations under standard Lipschitz smoothness. If second-order Lipschitz smoothness is added, then one can use momentum to non-trivially improve the complexity $T = O(n\varepsilon^{-1.5} + n^{0.75}\varepsilon^{-1.75})$ [1].

In the online stochastic setting, SGD finds a point $x$ with $\|\nabla f(x)\| \le \varepsilon$ in $T = O(\varepsilon^{-4})$ stochastic gradient iterations under the standard Lipschitz smoothness assumption. Perhaps surprisingly, without using momentum, one can already improve this rate to $T = O(\varepsilon^{-3.5})$ (using SGD) [4], $T = O(\varepsilon^{-3.333})$ [34], or even to $T = O(\varepsilon^{-3.25})$ [2] if second-order Lipschitz smoothness is added. It is unclear whether such rates can be improved using momentum. We stress that, even if "improved rates" can be obtained using momentum, one also needs to prove from the lower-bound side that such "improved rates" cannot be obtained by any momentum-free method.

[20] When $\sigma$ is large, even if $n$ is large, the iteration complexity of [28, 32] becomes $O(\sigma/\varepsilon^2)$. In this regime, almost all variance-reduction methods, including SVRG and `Katyusha`, can be shown to satisfy $\varepsilon \le O(\frac{\sqrt{\sigma}}{\sqrt{T}})$ within the first epoch, if the learning rates are appropriately chosen. Therefore, `Katyusha` and SVRG are no slower than [28, 32].

[21] Nitanda's method is usually not considered as an accelerated method, since it requires mini-batch size to be very large in order to be accelerated. If mini-batch is large then one can use full-gradient method directly and acceleration is trivial. This is confirmed by [31, Section IV.F]. In contrast, our acceleration holds even if mini-batch size is 1.

# B Experiment Details

The datasets we used in this paper are downloaded from the LibSVM website [23]:

- the adult (a9a) dataset (32,561 samples and 123 features).
- the web (w8a) dataset (49,749 samples and 300 features).
- the covtype (binary.scale) dataset (581,012 samples and 54 features).
- the mnist (class 1) dataset (60,000 samples and 780 features).
- the rcv1 (train.binary) dataset (20,242 samples and 47,236 features).
- the sensit (combined) dataset (78,823 samples and 100 features).

To make easier comparison across datasets, we scale every vector by the average Euclidean norm of all the vectors in the dataset. In other words, we ensure that the data vectors have an average Euclidean norm 1. This step is for comparison only and not necessary in practice.

**Parameter-tuning details.** We select learning rates from the set $\{10^{-k}, 2\times 10^{-k}, 5\times 10^{-k} : k\in\mathbb{Z}\}$, and select regularizer weights (for APCG) from the set $\{10^{-k} : k\in\mathbb{Z}\}$. We have fully automated the parameter tuning procedure to ensure a fair and strong comparison.

While the learning rates were explicitly defined for SVRG and SAGA, there were implicit for all accelerated methods. For Catalyst, the learning rate is in fact their $\alpha_0$ in the paper [35]. Instead of choosing it to be the theory-predicted value, we multiply it with an extra factor to be tuned and call this factor the "learning rate". Similarly, for `Katyusha` and `Katyusha`$^{\mathsf{ns}}$, we multiply the theory-predicted $\tau_1$ with an extra factor and this serves as a learning rate. For APCG, we use their Algorithm 1 in the paper and multiply their theory-predicted $\mu$ with an extra factor.

For Catalyst, in principle one also has to tune the stopping criterion. After communicating with an author of Catalyst, we learned that one can terminate the inner loop whenever the duality gap becomes no more than, say one fourth, of the last duality gap from the previous epoch [35]. This stopping criterion was also found by the authors of [6] to be a good choice for reduction-based methods.

# C Appendix for Section 4

## C.1 Proof of Theorem 4.1

*Proof of Theorem 4.1.* First of all, the parameter choices satisfy the presumptions in Lemma 2.6, so again by defining $D_k \stackrel{\text{def}}{=} F(y_k) - F(x^*)$ and $\widetilde{D}^s \stackrel{\text{def}}{=} F(\widetilde{x}^s) - F(x^*)$, we can rewrite Lemma 2.7 as follows:

$$0 \le \frac{\alpha_s(1-\tau_{1,s}-\tau_2)}{\tau_{1,s}} D_k - \frac{\alpha_s}{\tau_{1,s}}\mathbb{E}\big[D_{k+1}\big] + \frac{\alpha_s\tau_2}{\tau_{1,s}}\widetilde{D}^s + \frac{1}{2}\|z_k - x^*\|^2 - \frac{1}{2}\mathbb{E}\big[\|z_{k+1}-x^*\|^2\big] \enspace .$$

Summing up the above inequality for all the iterations $k = sm, sm+1, \dots, sm+m-1$, we have

$$\mathbb{E}\Big[\alpha_s\frac{1-\tau_{1,s}-\tau_2}{\tau_{1,s}} D_{(s+1)m} + \alpha_s\frac{\tau_{1,s}+\tau_2}{\tau_{1,s}}\sum_{j=1}^m D_{sm+j}\Big]$$
$$\le \alpha_s\frac{1-\tau_{1,s}-\tau_2}{\tau_{1,s}} D_{sm} + \alpha_s\frac{\tau_2}{\tau_{1,s}} m\widetilde{D}^s + \frac{1}{2}\|z_{sm}-x^*\|^2 - \frac{1}{2}\mathbb{E}\big[\|z_{(s+1)m}-x^*\|^2\big] \enspace . \tag{C.1}$$

Note that in the above inequality we have assumed all the randomness in the first $s-1$ epochs are fixed and the only source of randomness comes from epoch $s$.

If we define $\widetilde{x}^s = \frac{1}{m}\sum_{j=1}^m y_{(s-1)m+j}$, then by the convexity of function $F(\cdot)$ we have $m\widetilde{D}^s \le \sum_{j=1}^n D_{(s-1)m+j}$. Therefore, using the parameter choice $\alpha_s = \frac{1}{3\tau_{1,s}L}$, for every $s\ge 1$ we can derive

from (C.1) that

$$\mathbb{E}\Big[\frac{1}{\tau_{1,s}^2}D_{(s+1)m} + \frac{\tau_{1,s}+\tau_2}{\tau_{1,s}^2}\sum_{j=1}^{m-1} D_{sm+j}\Big]$$

$$\leq \frac{1-\tau_{1,s}}{\tau_{1,s}^2}D_{sm} + \frac{\tau_2}{\tau_{1,s}^2}\sum_{j=1}^{m-1} D_{(s-1)m+j} + \frac{3L}{2}\|z_{sm}-x^*\|^2 - \frac{3L}{2}\mathbb{E}\big[\|z_{(s+1)m}-x^*\|^2\big] \ . \tag{C.2}$$

For the base case $s=0$, we can also rewrite (C.1) as

$$\mathbb{E}\Big[\frac{1}{\tau_{1,0}^2}D_{m} + \frac{\tau_{1,0}+\tau_2}{\tau_{1,0}^2}\sum_{j=1}^{m-1} D_{j}\Big]$$

$$\leq \frac{1-\tau_{1,0}-\tau_2}{\tau_{1,0}^2}D_0 + \frac{\tau_2 m}{\tau_{1,0}^2}\widetilde{D}^0 + \frac{3L}{2}\|z_0-x^*\|^2 - \frac{3L}{2}\mathbb{E}\big[\|z_m-x^*\|^2\big] \ . \tag{C.3}$$

At this point, if we choose $\tau_{1,s}=\frac{2}{s+4}\leq\frac{1}{2}$, it satisfies

$$\frac{1}{\tau_{1,s}^2} \geq \frac{1-\tau_{1,s+1}}{\tau_{1,s+1}^2} \quad\text{and}\quad \frac{\tau_{1,s}+\tau_2}{\tau_{1,s}^2} \geq \frac{\tau_2}{\tau_{1,s+1}^2} \ .$$

Using these two inequalities, we can telescope (C.3) and (C.2) for all $s=0,1,\dots,S-1$. We obtain in the end that

$$\mathbb{E}\Big[\frac{1}{\tau_{1,S-1}^2}D_{Sm} + \frac{\tau_{1,S-1}+\tau_2}{\tau_{1,S-1}^2}\sum_{j=1}^{m-1} D_{(S-1)m+j} + \frac{3L}{2}\|z_{Sm}-z^*\|^2\Big]$$

$$\leq \frac{1-\tau_{1,0}-\tau_2}{\tau_{1,0}^2}D_0 + \frac{\tau_2 m}{\tau_{1,0}^2}\widetilde{D}^0 + \frac{3L}{2}\|z_0-x^*\|^2 \tag{C.4}$$

Since we have $\widetilde{D}^S \leq \frac{1}{m}\sum_{j=1}^m D_{(S-1)m+j}$ which is no greater than $\frac{2\tau_{1,S-1}^2}{m}$ times the left hand side of (C.4), we conclude that

$$\mathbb{E}\big[F(\widetilde{x}^S)-F(x^*)\big] = \mathbb{E}\big[\widetilde{D}^S\big] \leq O\big(\frac{\tau_{1,S}^2}{m}\big)\cdot\Big(\frac{1-\tau_{1,0}-\tau_2}{\tau_{1,0}^2}D_0 + \frac{\tau_2 m}{\tau_{1,0}^2}\widetilde{D}^0 + \frac{3L}{2}\|z_0-x^*\|^2\Big)$$

$$= O\big(\frac{1}{mS^2}\big)\cdot\Big(m\big(F(x_0)-F(x^*)\big) + L\|x_0-x^*\|^2\Big) \ . \qquad \square$$

# D Appendix for Section 5

## D.1 One-Iteration Analysis

Similar as Section 2.1, we first analyze the behavior of `Katyusha1` in a single iteration (i.e., for a fixed $k$). We view $y_k, z_k$ and $x_{k+1}$ as fixed in this section so the only randomness comes from the choice of $i$ in iteration $k$. We abbreviate in this subsection by $\widetilde{x}=\widetilde{x}^s$ where $s$ is the epoch that iteration $k$ belongs to, and denote by $\sigma_{k+1}^2 \stackrel{\text{def}}{=} \|\nabla f(x_{k+1}) - \widetilde{\nabla}_{k+1}\|^2$.

Our first lemma is analogous to Lemma 2.3, where note that we have replaced the use of $L$ in Lemma 2.3 with $L_\diamond \geq L$:

**Lemma D.1** (proximal gradient descent). *If* $L_\diamond \geq L$ *and*

$$y_{k+1} = \arg\min_y \Big\{ \frac{3L_\diamond}{2}\|y - x_{k+1}\|^2 + \langle \widetilde{\nabla}_{k+1}, y - x_{k+1}\rangle + \psi(y) - \psi(x_{k+1}) \Big\} \ , \qquad \textit{and}$$

$$\mathsf{Prog}(x_{k+1}) \stackrel{\text{def}}{=} -\min_y \Big\{ \frac{3L_\diamond}{2}\|y - x_{k+1}\|^2 + \langle \widetilde{\nabla}_{k+1}, y - x_{k+1}\rangle + \psi(y) - \psi(x_{k+1}) \Big\} \geq 0 \ ,$$

*we have (where the expectation is only over the randomness of* $\widetilde{\nabla}_{k+1}$*)*

$$F(x_{k+1}) - \mathbb{E}\big[F(y_{k+1})\big] \geq \mathbb{E}\big[\mathsf{Prog}(x_{k+1})\big] - \frac{1}{4L_\diamond}\mathbb{E}\big[\sigma_{k+1}^2\big] \ .$$

*Proof.*

$$\begin{aligned}
\mathsf{Prog}(x_{k+1}) &= -\min_y \Big\{ \frac{3L_\diamond}{2}\|y - x_{k+1}\|^2 + \langle \widetilde{\nabla}_{k+1}, y - x_{k+1}\rangle + \psi(y) - \psi(x_{k+1}) \Big\} \\
&\stackrel{①}{=} -\Big( \frac{3L_\diamond}{2}\|y_{k+1} - x_{k+1}\|^2 + \langle \widetilde{\nabla}_{k+1}, y_{k+1} - x_{k+1}\rangle + \psi(y_{k+1}) - \psi(x_{k+1}) \Big) \\
&= -\Big( \frac{L_\diamond}{2}\|y_{k+1} - x_{k+1}\|^2 + \langle \nabla f(x_{k+1}), y_{k+1} - x_{k+1}\rangle + \psi(y_{k+1}) - \psi(x_{k+1}) \Big) \\
&\quad + \Big( \langle \nabla f(x_{k+1}) - \widetilde{\nabla}_{k+1}, y_{k+1} - x_{k+1}\rangle - L_\diamond\|y_{k+1} - x_{k+1}\|^2 \Big) \\
&\stackrel{②}{\leq} -\Big( f(y_{k+1}) - f(x_{k+1}) + \psi(y_{k+1}) - \psi(x_{k+1}) \Big) + \frac{1}{4L_\diamond}\|\nabla f(x_{k+1}) - \widetilde{\nabla}_{k+1}\|^2 \ .
\end{aligned}$$

Above, ① is by the definition of $y_{k+1}$, and ② uses the smoothness of function $f(\cdot)$, as well as Young's inequality $\langle a, b\rangle - \frac{1}{2}\|b\|^2 \leq \frac{1}{2}\|a\|^2$. Taking expectation on both sides we arrive at the desired result. □

The following lemma is analogous to Lemma 2.4. The main difference is that since we have not chosen a mini-batch of size $b$, one should expect the variance to decrease by a factor of $b$. Also, since we are in the non-uniform case one should expect the use of $L$ in Lemma 2.4 to be replaced with $\overline{L}$:

**Lemma D.2** (variance upper bound).

$$\mathbb{E}\big[\|\widetilde{\nabla}_{k+1} - \nabla f(x_{k+1})\|^2\big] \leq \frac{2\overline{L}}{b} \cdot \big( f(\widetilde{x}) - f(x_{k+1}) - \langle \nabla f(x_{k+1}), \widetilde{x} - x_{k+1}\rangle \big) \ .$$

*Proof.* Each $f_i(x)$, being convex and $L_i$-smooth, implies the following inequality which is classical in convex optimization and can be found for instance in Theorem 2.1.5 of the textbook of Nesterov [43].

$$\|\nabla f_i(x_{k+1}) - \nabla f_i(\widetilde{x})\|^2 \leq 2L_i \cdot \big( f_i(\widetilde{x}) - f_i(x_{k+1}) - \langle \nabla f_i(x_{k+1}), \widetilde{x} - x_{k+1}\rangle \big)$$

Therefore, taking expectation over the random choice of $i$, we have

$$\begin{aligned}
&\mathbb{E}\big[\|\widetilde{\nabla}_{k+1} - \nabla f(x_{k+1})\|^2\big] \\
&= \mathbb{E}_{S_k}\Big[\Big\|\Big(\frac{1}{b}\sum_{i\in S_k}\Big(\nabla f(\widetilde{x}) + \frac{1}{np_i}\big(\nabla f_i(x_{k+1}) - \nabla f_i(\widetilde{x})\big)\Big)\Big) - \nabla f(x_{k+1})\Big\|^2\Big] \\
&= \frac{1}{b}\mathbb{E}_{i\sim\mathcal{D}}\Big[\Big\|\Big(\nabla f(\widetilde{x}) + \frac{1}{np_i}\big(\nabla f_i(x_{k+1}) - \nabla f_i(\widetilde{x})\big)\Big) - \nabla f(x_{k+1})\Big\|^2\Big] \\
&= \frac{1}{b}\mathbb{E}_{i\sim\mathcal{D}}\Big[\Big\|\frac{1}{np_i}\big(\nabla f_i(x_{k+1}) - \nabla f_i(\widetilde{x})\big) - \big(\nabla f(x_{k+1}) - f(\widetilde{x})\big)\Big\|^2\Big] \\
&\overset{①}{\le} \frac{1}{b}\mathbb{E}_{i\sim\mathcal{D}}\Big[\Big\|\frac{1}{np_i}\big(\nabla f_i(x_{k+1}) - \nabla f_i(\widetilde{x})\big)\Big\|^2\Big] \\
&\overset{②}{\le} \frac{1}{b}\cdot\sum_{i\in[n]}\frac{2L_i}{n^2p_i}\Big(f_i(\widetilde{x}) - f_i(x_{k+1}) - \langle\nabla f_i(x_{k+1}), \widetilde{x} - x_{k+1}\rangle\Big) \\
&= \frac{2\overline{L}}{b}\cdot\big(f(\widetilde{x}) - f(x_{k+1}) - \langle\nabla f(x_{k+1}), \widetilde{x} - x_{k+1}\rangle\big) \ .
\end{aligned}$$

Above, ① is because for any random vector $\zeta \in \mathbb{R}^d$, it holds that $\mathbb{E}\|\zeta - \mathbb{E}\zeta\|^2 = \mathbb{E}\|\zeta\|^2 - \|\mathbb{E}\zeta\|^2$; ② follows from the first inequality in this proof. □

The next lemma is completely identical to Lemma 2.5 so we skip the proof.

**Lemma D.3** (proximal mirror descent)**.** *Suppose $\psi(\cdot)$ is $\sigma$-SC. Then, fixing $\widetilde{\nabla}_{k+1}$ and letting*

$$z_{k+1} = \arg\min_z\big\{\frac{1}{2}\|z - z_k\|^2 + \alpha\langle\widetilde{\nabla}_{k+1}, z - z_k\rangle + \alpha\psi(z) - \alpha\psi(z_k)\big\} \ ,$$

*it satisfies for all $u \in \mathbb{R}^d$,*

$$\alpha\langle\widetilde{\nabla}_{k+1}, z_{k+1} - u\rangle + \alpha\psi(z_{k+1}) - \alpha\psi(u) \le -\frac{1}{2}\|z_k - z_{k+1}\|^2 + \frac{1}{2}\|z_k - u\|^2 - \frac{1+\alpha\sigma}{2}\|z_{k+1} - u\|^2 \ .$$

The following lemma combines Lemma D.1, Lemma D.2 and Lemma D.3 all together, using the special choice of $x_{k+1}$ which is a convex combination of $y_k, z_k$ and $\widetilde{x}$:

**Lemma D.4** (coupling step 1)**.** *If $x_{k+1} = \tau_1 z_k + \tau_2\widetilde{x} + (1 - \tau_1 - \tau_2)y_k$, where $\tau_1 \le \frac{1}{3\alpha L_\diamond}$ and $\tau_2 = \frac{\overline{L}}{2L_\diamond b}$,*

$$\begin{aligned}
&\alpha\langle\nabla f(x_{k+1}), z_k - u\rangle - \alpha\psi(u) \\
&\le \frac{\alpha}{\tau_1}\Big(F(x_{k+1}) - \mathbb{E}\big[F(y_{k+1})\big] + \tau_2 F(\widetilde{x}) - \tau_2 f(x_{k+1}) - \tau_2\langle\nabla f(x_{k+1}), \widetilde{x} - x_{k+1}\rangle\Big) \\
&\quad + \frac{1}{2}\|z_k - u\|^2 - \frac{1+\alpha\sigma}{2}\mathbb{E}\big[\|z_{k+1} - u\|^2\big] + \frac{\alpha(1-\tau_1-\tau_2)}{\tau_1}\psi(y_k) - \frac{\alpha}{\tau_1}\psi(x_{k+1}) \ .
\end{aligned}$$

*Proof.* We first apply Lemma D.3 and get

$$\begin{aligned}
&\alpha\langle\widetilde{\nabla}_{k+1}, z_k - u\rangle + \alpha\psi(z_{k+1}) - \alpha\psi(u) \\
&= \alpha\langle\widetilde{\nabla}_{k+1}, z_k - z_{k+1}\rangle + \alpha\langle\widetilde{\nabla}_{k+1}, z_{k+1} - u\rangle + \alpha\psi(z_{k+1}) - \alpha\psi(u) \\
&\le \alpha\langle\widetilde{\nabla}_{k+1}, z_k - z_{k+1}\rangle - \frac{1}{2}\|z_k - z_{k+1}\|^2 + \frac{1}{2}\|z_k - u\|^2 - \frac{1+\alpha\sigma}{2}\|z_{k+1} - u\|^2 \ . \qquad \text{(D.1)}
\end{aligned}$$

By defining $v \stackrel{\text{def}}{=} \tau_1 z_{k+1} + \tau_2 \widetilde{x} + (1-\tau_1-\tau_2) y_k$, we have $x_{k+1} - v = \tau_1(z_k - z_{k+1})$ and therefore

$$\begin{aligned}
&\mathbb{E}\Big[\alpha\langle\widetilde{\nabla}_{k+1}, z_k - z_{k+1}\rangle - \frac{1}{2}\|z_k - z_{k+1}\|^2\Big] = \mathbb{E}\Big[\frac{\alpha}{\tau_1}\langle\widetilde{\nabla}_{k+1}, x_{k+1} - v\rangle - \frac{1}{2\tau_1^2}\|x_{k+1} - v\|^2\Big] \\
&= \mathbb{E}\Big[\frac{\alpha}{\tau_1}\Big(\langle\widetilde{\nabla}_{k+1}, x_{k+1} - v\rangle - \frac{1}{2\alpha\tau_1}\|x_{k+1} - v\|^2 - \psi(v) + \psi(x_{k+1})\Big) + \frac{\alpha}{\tau_1}\Big(\psi(v) - \psi(x_{k+1})\Big)\Big] \\
&\overset{①}{\leq} \mathbb{E}\Big[\frac{\alpha}{\tau_1}\Big(\langle\widetilde{\nabla}_{k+1}, x_{k+1} - v\rangle - \frac{3L_\diamond}{2}\|x_{k+1} - v\|^2 - \psi(v) + \psi(x_{k+1})\Big) + \frac{\alpha}{\tau_1}\Big(\psi(v) - \psi(x_{k+1})\Big)\Big] \\
&\overset{②}{\leq} \mathbb{E}\Big[\frac{\alpha}{\tau_1}\Big(F(x_{k+1}) - F(y_{k+1}) + \frac{1}{4L_\diamond}\sigma_{k+1}^2\Big) + \frac{\alpha}{\tau_1}\Big(\psi(v) - \psi(x_{k+1})\Big)\Big] \\
&\overset{③}{\leq} \mathbb{E}\Big[\frac{\alpha}{\tau_1}\Big(F(x_{k+1}) - F(y_{k+1}) + \frac{\overline{L}}{2L_\diamond b}\big(f(\widetilde{x}) - f(x_{k+1}) - \langle\nabla f(x_{k+1}), \widetilde{x} - x_{k+1}\rangle\big)\Big) \\
&\qquad + \frac{\alpha}{\tau_1}\Big(\tau_1\psi(z_{k+1}) + \tau_2\psi(\widetilde{x}) + (1-\tau_1-\tau_2)\psi(y_k) - \psi(x_{k+1})\Big)\Big] \ . \qquad \text{(D.2)}
\end{aligned}$$

Above, ① uses our choice $\tau_1 \leq \frac{1}{3\alpha L}$, ② uses Lemma D.1, ③ uses Lemma D.2 together with the convexity of $\psi(\cdot)$ and the definition of $v$. Finally, noticing that $\mathbb{E}[\langle\widetilde{\nabla}_{k+1}, z_k - u\rangle] = \langle\nabla f(x_{k+1}), z_k - u\rangle$ and $\tau_2 = \frac{1}{2}$, we obtain the desired inequality by combining (D.1) and (D.2). □

The next lemma simplifies the left hand side of Lemma D.4 using the convexity of $f(\cdot)$, and gives an inequality that relates the objective-distance-to-minimizer quantities $F(y_k) - F(x^*)$, $F(y_{k+1}) - F(x^*)$, and $F(\widetilde{x}) - F(x^*)$ to the point-distance-to-minimizer quantities $\|z_k - x^*\|^2$ and $\|z_{k+1} - x^*\|^2$.

**Lemma D.5** (coupling step 2)**.** *Under the same choices of $\tau_1, \tau_2$ as in Lemma D.4, we have*

$$\begin{aligned}
0 \leq \frac{\alpha(1-\tau_1-\tau_2)}{\tau_1}(F(y_k) - F(x^*)) - \frac{\alpha}{\tau_1}\big(\mathbb{E}\big[F(y_{k+1})\big] - F(x^*)\big) + \frac{\alpha\tau_2}{\tau_1}\big(F(\widetilde{x}) - F(x^*)\big) \\
+ \frac{1}{2}\|z_k - x^*\|^2 - \frac{1+\alpha\sigma}{2}\mathbb{E}\big[\|z_{k+1} - x^*\|^2\big] \ .
\end{aligned}$$

*Proof.* We first compute that

$$\begin{aligned}
&\alpha\big(f(x_{k+1}) - f(u)\big) \overset{①}{\leq} \alpha\langle\nabla f(x_{k+1}), x_{k+1} - u\rangle \\
&= \alpha\langle\nabla f(x_{k+1}), x_{k+1} - z_k\rangle + \alpha\langle\nabla f(x_{k+1}), z_k - u\rangle \\
&\overset{②}{=} \frac{\alpha\tau_2}{\tau_1}\langle\nabla f(x_{k+1}), \widetilde{x} - x_{k+1}\rangle + \frac{\alpha(1-\tau_1-\tau_2)}{\tau_1}\langle\nabla f(x_{k+1}), y_k - x_{k+1}\rangle + \alpha\langle\nabla f(x_{k+1}), z_k - u\rangle \\
&\overset{③}{\leq} \frac{\alpha\tau_2}{\tau_1}\langle\nabla f(x_{k+1}), \widetilde{x} - x_{k+1}\rangle + \frac{\alpha(1-\tau_1-\tau_2)}{\tau_1}(f(y_k) - f(x_{k+1})) + \alpha\langle\nabla f(x_{k+1}), z_k - u\rangle \ .
\end{aligned}$$

Above, ① uses the convexity of $f(\cdot)$, ② uses the choice that $x_{k+1} = \tau_1 z_k + \tau_2\widetilde{x} + (1-\tau_1-\tau_2)y_k$, and ③ uses the convexity of $f(\cdot)$ again. By applying Lemma D.4 to the above inequality, we have

$$\begin{aligned}
&\alpha\big(f(x_{k+1}) - F(u)\big) \leq \frac{\alpha(1-\tau_1-\tau_2)}{\tau_1}(F(y_k) - f(x_{k+1})) \\
&+ \frac{\alpha}{\tau_1}\Big(F(x_{k+1}) - \mathbb{E}\big[F(y_{k+1})\big] + \tau_2 F(\widetilde{x}) - \tau_2 f(x_{k+1})\Big) + \frac{1}{2}\|z_k - u\|^2 - \frac{1+\alpha\sigma}{2}\mathbb{E}\big[\|z_{k+1} - u\|^2\big] - \frac{\alpha}{\tau_1}\psi(x_{k+1})
\end{aligned}$$

which implies

$$\begin{aligned}
&\alpha\big(F(x_{k+1}) - F(u)\big) \leq \frac{\alpha(1-\tau_1-\tau_2)}{\tau_1}(F(y_k) - F(x_{k+1})) \\
&\quad + \frac{\alpha}{\tau_1}\Big(F(x_{k+1}) - \mathbb{E}\big[F(y_{k+1})\big] + \tau_2 F(\widetilde{x}) - \tau_2 F(x_{k+1})\Big) + \frac{1}{2}\|z_k - u\|^2 - \frac{1+\alpha\sigma}{2}\mathbb{E}\big[\|z_{k+1} - u\|^2\big] \ .
\end{aligned}$$

After rearranging and setting $u = x^*$, the above inequality yields

$$0 \le \frac{\alpha(1-\tau_1-\tau_2)}{\tau_1}\big(F(y_k)-F(x^*)\big) - \frac{\alpha}{\tau_1}\big(\mathbb{E}\big[F(y_{k+1})-F(x^*)\big]\big) + \frac{\alpha\tau_2}{\tau_1}\big(F(\widetilde{x})-F(x^*)\big) + \frac{1}{2}\|z_k - x^*\|^2 - \frac{1+\alpha\sigma}{2}\mathbb{E}\big[\|z_{k+1}-x^*\|^2\big] \ . \qquad \square$$

## D.2 Proof of Theorem 5.4

We are now ready to combine the analyses across iterations, and derive our final Theorem 5.4. Our proof next requires a careful telescoping of Lemma D.5 together with our specific parameter choices.

*Proof of Theorem 5.4.* Define $D_k \stackrel{\text{def}}{=} F(y_k) - F(x^*)$, $\widetilde{D}^s \stackrel{\text{def}}{=} F(\widetilde{x}^s) - F(x^*)$, and rewrite Lemma D.5:

$$0 \le \frac{(1-\tau_1-\tau_2)}{\tau_1}D_k - \frac{1}{\tau_1}D_{k+1} + \frac{\tau_2}{\tau_1}\mathbb{E}\big[\widetilde{D}^s\big] + \frac{1}{2\alpha}\|z_k - x^*\|^2 - \frac{1+\alpha\sigma}{2\alpha}\mathbb{E}\big[\|z_{k+1}-x^*\|^2\big] \ .$$

At this point, let us $\theta$ be an arbitrary value in $\big[1, 1+\alpha\sigma\big]$ and multiply the above inequality by $\theta^j$ for each $k = sm + j$. Then, we sum up the resulting $m$ inequalities for all $j = 0, 1, \dots, m-1$:

$$0 \le \mathbb{E}\Big[\frac{(1-\tau_1-\tau_2)}{\tau_1}\sum_{j=0}^{m-1} D_{sm+j}\cdot\theta^j - \frac{1}{\tau_1}\sum_{j=0}^{m-1} D_{sm+j+1}\cdot\theta^j\Big] + \frac{\tau_2}{\tau_1}\widetilde{D}^s \cdot \sum_{j=0}^{m-1}\theta^j + \frac{1}{2\alpha}\|z_{sm} - x^*\|^2 - \frac{\theta^m}{2\alpha}\big[\|z_{(s+1)m} - x^*\|^2\big] \ .$$

Note that in the above inequality we have assumed all the randomness in the first $s-1$ epochs are fixed and the only source of randomness comes from epoch $s$. We can rearrange the terms in the above inequality and get

$$\mathbb{E}\Big[\frac{\tau_1+\tau_2-(1-1/\theta)}{\tau_1}\sum_{j=1}^{m} D_{sm+j}\cdot\theta^j\Big] \le \frac{(1-\tau_1-\tau_2)}{\tau_1}\Big(D_{sm} - \theta^m\mathbb{E}\big[D_{(s+1)m}\big]\Big) + \frac{\tau_2}{\tau_1}\widetilde{D}^s\cdot\sum_{j=0}^{m-1}\theta^j + \frac{1}{2\alpha}\|z_{sm}-x^*\|^2 - \frac{\theta^m}{2\alpha}\mathbb{E}\big[\|z_{(s+1)m}-x^*\|^2\big] \ .$$

Using the special choice that $\widetilde{x}^{s+1} = \big(\sum_{j=0}^{m-1}\theta^j\big)^{-1}\cdot\sum_{j=0}^{m-1} y_{sm+j+1}\cdot\theta^j$ and the convexity of $F(\cdot)$, we derive that $\widetilde{D}^{s+1} \le \big(\sum_{j=0}^{m-1}\theta^j\big)^{-1}\cdot\sum_{j=0}^{m-1} D_{sm+j+1}\cdot\theta^j$. Substituting this into the above inequality, we get

$$\frac{\tau_1+\tau_2-(1-1/\theta)}{\tau_1}\theta\mathbb{E}\big[\widetilde{D}^{s+1}\big]\cdot\sum_{j=0}^{m-1}\theta^j \le \frac{(1-\tau_1-\tau_2)}{\tau_1}\Big(D_{sm} - \theta^m\mathbb{E}\big[D_{(s+1)m}\big]\Big) + \frac{\tau_2}{\tau_1}\widetilde{D}^s\cdot\sum_{j=0}^{m-1}\theta^j + \frac{1}{2\alpha}\|z_{sm}-x^*\|^2 - \frac{\theta^m}{2\alpha}\mathbb{E}\big[\|z_{(s+1)m}-x^*\|^2\big] \ . \quad \text{(D.3)}$$

We consider two cases (and four subcases) next.

**Case 1.** Suppose $L \le \frac{\overline{L}m}{b}$. In this case, we choose

$$\tau_2 = \min\Big\{\frac{\overline{L}}{2Lb}, \frac{1}{2}\Big\} \in \Big[\frac{1}{2m}, \frac{1}{2}\Big] \quad \text{and} \quad L_\diamond = \frac{\overline{L}}{2b\tau_2} \ge L$$

**Case 1.1.** Suppose $\frac{m\sigma b}{\overline{L}} \le \frac{3}{8}$. In this subcase, we choose

$$\alpha = \frac{\sqrt{b}}{\sqrt{6m\sigma\overline{L}}} \ , \qquad \tau_1 = \frac{1}{3\alpha L_\diamond} = 4m\alpha\sigma\tau_2 = \frac{\sqrt{8\tau_2^2 bm\sigma}}{\sqrt{3\overline{L}}} \in [0, \tau_2] \subseteq [0, \frac{1}{2}] \ , \quad \text{and} \quad \theta = 1 + \alpha\sigma$$

We have

$$\alpha\sigma = \frac{1}{\sqrt{6m^2}} \frac{\sqrt{b\sigma m}}{\sqrt{\overline{L}}} \le \frac{1}{4m}$$

and therefore the following inequality holds:

$$\tau_2(\theta^{m-1} - 1) + (1 - 1/\theta) = \tau_2((1 + \alpha\sigma)^{m-1} - 1) + (1 - \frac{1}{1 + \alpha\sigma}) \le 2\tau_2 m\alpha\sigma + \alpha\sigma \le 4\tau_2 m\alpha\sigma = \tau_1 \ .$$

In other words, we have $\tau_1 + \tau_2 - (1 - 1/\theta) \ge \tau_2\theta^{m-1}$ and thus (D.3) implies that

$$\mathbb{E}\Big[\frac{\tau_2}{\tau_1}\widetilde{D}^{s+1} \cdot \sum_{j=0}^{m-1} \theta^j + \frac{1 - \tau_1 - \tau_2}{\tau_1} D_{(s+1)m} + \frac{1}{2\alpha}\|z_{(s+1)m} - x^*\|^2\Big]$$
$$\le \theta^{-m} \cdot \Big(\frac{\tau_2}{\tau_1}\widetilde{D}^{s} \cdot \sum_{j=0}^{m-1} \theta^j + \frac{1 - \tau_1 - \tau_2}{\tau_1} D_{sm} + \frac{1}{2\alpha}\|z_{sm} - x^*\|^2\Big) \ .$$

If we telescope the above inequality over all epochs $s = 0, 1, \dots, S - 1$, we obtain

$$\begin{aligned}
\mathbb{E}\big[F(x^{\mathsf{out}}) - F(x^*)\big] &\overset{①}{\le} \frac{1}{\tau_2 m + (1 - \tau_1 - \tau_2)}\mathbb{E}\big[\tau_2\widetilde{D}^S \cdot \sum_{j=0}^{m-1} \theta^j + \big(1 - \tau_1 - \tau_2\big)D_{Sm}\big] \\
&\overset{②}{\le} \theta^{-Sm} \cdot O\Big(\widetilde{D}^0 + D_0 + \frac{\tau_1}{\alpha\tau_2 m}\|x_0 - x^*\|^2\Big) \\
&\overset{③}{\le} \theta^{-Sm} \cdot O\Big(1 + \frac{\tau_1}{\alpha\tau_2 m\sigma}\Big) \cdot (F(x_0) - F(x^*)) \\
&\overset{④}{=} O((1 + \alpha\sigma)^{-Sm}) \cdot \big(F(x_0) - F(x^*)\big) \ . 
\end{aligned} \tag{D.4}$$

Above, inequality ① uses the choice $x^{\mathsf{out}} = \frac{\tau_2 m\widetilde{x}^S + (1-\tau_1-\tau_2)y_{Sm}}{\tau_2 m + (1-\tau_1-\tau_2)}$, the convexity of $F(\cdot)$, and the fact $\sum_{j=0}^{m-1} \theta^j \ge m$; inequality ② uses the fact that $\sum_{j=0}^{m-1} \theta^j \le O(m)$ (because $\alpha\sigma \le \frac{1}{4m}$), and the fact that $\tau_2 \ge \frac{1}{2m}$; inequality ③ uses the strong convexity of $F(\cdot)$ which implies $F(x_0) - F(x^*) \ge \frac{\sigma}{2}\|x_0 - x^*\|^2$; and inequality ④ uses our choice of $\tau_1$.

**Case 1.2.** Suppose $\frac{m\sigma b}{\overline{L}} > \frac{3}{8}$. In this case, we choose

$$\tau_1 = \tau_2 \quad \text{and} \quad \alpha = \frac{1}{3\tau_1 L_\diamond} = \frac{2b}{3\overline{L}} \ge \frac{1}{4\sigma m} \ , \quad \theta = 1 + \frac{1}{4m}$$

(Note that we can choose $\theta = 1 + \frac{1}{4m}$ because $\frac{1}{4m} \le \alpha\sigma$.)

Under these parameter choices, we can calculate that

$$\frac{(\tau_1 + \tau_2 - (1 - 1/\theta))\theta}{\tau_2} = 2 - \frac{1 - 2\tau_2}{4m\tau_2} \ge \frac{3}{2} > \frac{5}{4} \quad \text{and} \quad \theta^m \ge \frac{5}{4}$$

thus (D.3) implies that

$$\mathbb{E}\Big[\frac{\tau_2}{\tau_1}\widetilde{D}^{s+1} \cdot \sum_{j=0}^{m-1} \theta^j + \frac{1 - \tau_1 - \tau_2}{\tau_1} D_{(s+1)m} + \frac{1}{2\alpha}\|z_{(s+1)m} - x^*\|^2\Big]$$
$$\le \frac{4}{5} \cdot \Big(\frac{\tau_2}{\tau_1}\widetilde{D}^{s} \cdot \sum_{j=0}^{m-1} \theta^j + \frac{1 - \tau_1 - \tau_2}{\tau_1} D_{sm} + \frac{1}{2\alpha}\|z_{sm} - x^*\|^2\Big) \ .$$

If we telescope the above inequality over all epochs $s = 0, 1, \ldots, S-1$, we obtain

$$\mathbb{E}\big[F(x^{\mathsf{out}}) - F(x^*)\big] \overset{①}{\le} \frac{1}{\tau_2 m + (1-\tau_1-\tau_2)} \mathbb{E}\big[\tau_2 \widetilde{D}^S \cdot \sum_{j=0}^{m-1} \theta^j + \big(1-\tau_1-\tau_2\big) D_{Sm}\big]$$

$$\overset{②}{\le} \big(\frac{5}{4}\big)^{-S} \cdot O\Big(\widetilde{D}^0 + D_0 + \frac{\tau_1}{\alpha \tau_2 m}\|x_0 - x^*\|^2\Big)$$

$$\overset{③}{\le} \big(\frac{5}{4}\big)^{-S} \cdot O\Big(1 + \frac{\tau_1}{\alpha \tau_2 m \sigma}\Big) \cdot (F(x_0) - F(x^*))$$

$$\overset{④}{=} O((5/4)^{-S}) \cdot \big(F(x_0) - F(x^*)\big) \ . \tag{D.5}$$

Above, inequality ① uses the choice $x^{\mathsf{out}} = \frac{\tau_2 m \widetilde{x}^S + (1-\tau_1-\tau_2) y_{Sm}}{\tau_2 m + (1-\tau_1-\tau_2)}$, the convexity of $F(\cdot)$, and the fact $\sum_{j=0}^{m-1} \theta^j \ge m$; inequality ② uses the fact that $\sum_{j=0}^{m-1} \theta^j \le O(m)$, and the fact that $\tau_2 \ge \frac{1}{2m}$; inequality ③ uses the strong convexity of $F(\cdot)$ which implies $F(x_0) - F(x^*) \ge \frac{\sigma}{2}\|x_0 - x^*\|^2$; and inequality ④ uses our choice of $\tau_1$ and $\alpha$.

**Case 2.** Suppose $L > \frac{\overline{L} m}{b}$. In this case, we choose

$$L_\diamond = L \quad \text{and} \quad \tau_2 = \frac{\overline{L}}{2 L_\diamond b} = \frac{\overline{L}}{2 L b} \in \big[0, \frac{1}{2m}\big]$$

**Case 2.1.** Suppose $\frac{m^2 \sigma}{L} \le \frac{3}{8}$. In this subcase, we choose

$$\alpha = \frac{1}{\sqrt{6\sigma L}} \ , \qquad \tau_1 = \frac{1}{3\alpha L} = 2\alpha\sigma = \frac{\sqrt{2\sigma}}{\sqrt{3L}} \in \big[0, \frac{1}{2m}\big] \ , \quad \theta = 1 + \alpha\sigma$$

We have $\alpha\sigma \le \frac{1}{4m}$ and therefore the following inequality holds:

$$\tau_2(\theta^{m-1} - 1) + (1 - 1/\theta) = \tau_2((1+\alpha\sigma)^{m-1} - 1) + (1 - \frac{1}{1+\alpha\sigma}) \le 2\tau_2 m \alpha\sigma + \alpha\sigma \le 2\alpha\sigma = \tau_1 \ .$$

In other words, we have $\tau_1 + \tau_2 - (1 - 1/\theta) \ge \tau_2 \theta^{m-1}$ and thus (D.3) implies that

$$\mathbb{E}\Big[\frac{\tau_2}{\tau_1} \widetilde{D}^{s+1} \cdot \sum_{j=0}^{m-1} \theta^j + \frac{1-\tau_1-\tau_2}{\tau_1} D_{(s+1)m} + \frac{1}{2\alpha}\|z_{(s+1)m} - x^*\|^2\Big]$$

$$\le \theta^{-m} \cdot \Big(\frac{\tau_2}{\tau_1} \widetilde{D}^{s} \cdot \sum_{j=0}^{m-1} \theta^j + \frac{1-\tau_1-\tau_2}{\tau_1} D_{sm} + \frac{1}{2\alpha}\|z_{sm} - x^*\|^2\Big) \ .$$

If we telescope the above inequality over all epochs $s = 0, 1, \ldots, S-1$, we obtain

$$\mathbb{E}\big[F(x^{\mathsf{out}}) - F(x^*)\big] \overset{①}{\le} \frac{1}{\tau_2 m + (1-\tau_1-\tau_2)} \mathbb{E}\big[\tau_2 \widetilde{D}^S \cdot \sum_{j=0}^{m-1} \theta^j + \big(1-\tau_1-\tau_2\big) D_{Sm}\big]$$

$$\overset{②}{\le} \theta^{-Sm} \cdot O\Big(\widetilde{D}^0 + D_0 + \frac{\tau_1}{\alpha}\|x_0 - x^*\|^2\Big)$$

$$\overset{③}{\le} \theta^{-Sm} \cdot O\Big(1 + \frac{\tau_1}{\alpha\sigma}\Big) \cdot (F(x_0) - F(x^*))$$

$$\overset{④}{=} O((1+\alpha\sigma)^{-Sm}) \cdot \big(F(x_0) - F(x^*)\big) \ . \tag{D.6}$$

Above, inequality ① uses the choice $x^{\mathsf{out}} = \frac{\tau_2 m \widetilde{x}^S + (1-\tau_1-\tau_2) y_{Sm}}{\tau_2 m + (1-\tau_1-\tau_2)}$, the convexity of $F(\cdot)$, and the fact $\sum_{j=0}^{m-1} \theta^j \ge m$; inequality ② uses the fact that $\sum_{j=0}^{m-1} \theta^j \le O(m)$ (because $\alpha\sigma \le \frac{1}{4m}$), and the fact that $\tau_2 m + (1-\tau_1-\tau_2) \geq 1-\tau_1+(m-1)\tau_2 \geq 1/2$; inequality ③ uses the strong convexity of $F(\cdot)$ which implies $F(x_0) - F(x^*) \geq \frac{\sigma}{2}\|x_0 - x^*\|^2$; and inequality ④ uses our choice of $\tau_1$.

**Case 2.2.** Suppose $\frac{m^2\sigma}{L} > \frac{3}{8}$. In this case, we choose

$$\tau_1 = \frac{1}{2m} \quad \text{and} \quad \alpha = \frac{1}{3\tau_1 L} = \frac{2m}{3L} > \frac{1}{4\sigma m} \enspace, \quad \theta = 1 + \frac{1}{4m}$$

(Note that we can choose $\theta = 1 + \frac{1}{4m}$ because $\frac{1}{4m} \leq \alpha\sigma$.)

Under these parameter choices, we can calculate that

$$\frac{(\tau_1+\tau_2-(1-1/\theta))\theta}{\tau_2} = \frac{\tau_1+\tau_2}{\tau_2} - \frac{1-2\tau_2}{4m\tau_2} \geq 1 + \frac{\tau_1 - 1/4m}{\tau_2} \geq \frac{3}{2} > \frac{5}{4} \quad \text{and} \quad \theta^m \geq \frac{5}{4}$$

thus (D.3) implies that

$$\begin{aligned}
&\mathbb{E}\Big[\frac{\tau_2}{\tau_1}\widetilde{D}^{s+1} \cdot \sum_{j=0}^{m-1}\theta^j + \frac{1-\tau_1-\tau_2}{\tau_1}D_{(s+1)m} + \frac{1}{2\alpha}\|z_{(s+1)m} - x^*\|^2\Big] \\
&\leq \frac{4}{5} \cdot \Big(\frac{\tau_2}{\tau_1}\widetilde{D}^{s} \cdot \sum_{j=0}^{m-1}\theta^j + \frac{1-\tau_1-\tau_2}{\tau_1}D_{sm} + \frac{1}{2\alpha}\|z_{sm} - x^*\|^2\Big) \enspace .
\end{aligned}$$

If we telescope the above inequality over all epochs $s = 0, 1, \dots, S-1$, we obtain

$$\begin{aligned}
\mathbb{E}\big[F(x^{\mathsf{out}}) - F(x^*)\big] &\overset{①}{\leq} \frac{1}{\tau_2 m + (1-\tau_1-\tau_2)}\mathbb{E}\big[\tau_2 \widetilde{D}^S \cdot \sum_{j=0}^{m-1}\theta^j + \big(1-\tau_1-\tau_2\big)D_{Sm}\big] \\
&\overset{②}{\leq} \big(\frac{5}{4}\big)^{-S} \cdot O\Big(\widetilde{D}^0 + D_0 + \frac{\tau_1}{\alpha}\|x_0 - x^*\|^2\Big) \\
&\overset{③}{\leq} \big(\frac{5}{4}\big)^{-S} \cdot O\Big(1 + \frac{\tau_1}{\alpha\sigma}\Big) \cdot (F(x_0) - F(x^*)) \\
&\overset{④}{=} O((5/4)^{-S}) \cdot \big(F(x_0) - F(x^*)\big) \enspace . \qquad \text{(D.7)}
\end{aligned}$$

Above, inequality ① uses the choice $x^{\mathsf{out}} = \frac{\tau_2 m \widetilde{x}^S + (1-\tau_1-\tau_2)y_{Sm}}{\tau_2 m + (1-\tau_1-\tau_2)}$, the convexity of $F(\cdot)$, and the fact $\sum_{j=0}^{m-1}\theta^j \geq m$; inequality ② uses the fact that $\sum_{j=0}^{m-1}\theta^j \leq O(m)$, and that $\tau_2 m + (1-\tau_1-\tau_2) \geq 1-\tau_1+(m-1)\tau_2 \geq 1/2$; inequality ③ uses the strong convexity of $F(\cdot)$ which implies $F(x_0) - F(x^*) \geq \frac{\sigma}{2}\|x_0 - x^*\|^2$; and inequality ④ uses our choice of $\tau_1$ and $\alpha$. □

# E Appendix for Section 6

In this section, we first include the complete pseudo-codes for `Katyusha2` and `Katyusha2`$^{\mathsf{ns}}$. Then, we provide a one-iteration analysis for both algorithms, in the same spirit as Section 2.1.

The final proofs of Theorem 6.1 and Theorem 6.2 are direct corollaries of such one-iteration analysis, where the details we have already given in Section 2.2 and in Section C.1 respectively.

## E.1 Pseudo-Codes

**Algorithm 4** `Katyusha2`$(x_0, S, \sigma, (L_1, \ldots, L_n))$

$m \leftarrow n$; $\overline{L} = (L_1 + \cdots + L_n)/n$;
$\tau_2 \leftarrow \frac{1}{2}$, $\tau_1 \leftarrow \min\big\{\sqrt{m\sigma/9\overline{L}}, \frac{1}{2}\big\}$, $\alpha \leftarrow \frac{1}{9\tau_1\overline{L}}$;
$y_0 = z_0 = \widetilde{x}^0 \leftarrow x_0$;
**for** $s \leftarrow 0$ **to** $S-1$ **do**
    $\mu^s \leftarrow \nabla f(\widetilde{x}^s)$;
    **for** $j \leftarrow 0$ **to** $m-1$ **do**
        $k \leftarrow (sm) + j$;
        $x_{k+1} \leftarrow \tau_1 z_k + \tau_2 \widetilde{x}^s + (1 - \tau_1 - \tau_2) y_k$;
        Pick $i$ randomly from $\{1, 2, \ldots, n\}$, each with probability $L_i/n\overline{L}$;
        $\widetilde{\nabla}_{k+1} \leftarrow \mu^s + \nabla f_i(x_{k+1}) - \nabla f_i(\widetilde{x}^s)$;
        $z_{k+1} = \arg\min_z \big\{\frac{1}{\alpha} V_{z_k}(z) + \langle \widetilde{\nabla}_{k+1}, z\rangle + \psi(z)\big\}$;
            ◇ *$V_x(y)$ is the Bregman divergence function, see Section 6*
        $y_{k+1} \leftarrow \arg\min_y \big\{\frac{9\overline{L}}{2}\|y - x_{k+1}\|^2 + \langle \widetilde{\nabla}_{k+1}, y\rangle + \psi(y)\big\}$;
    **end for**
    $\widetilde{x}^{s+1} \leftarrow \big(\sum_{j=0}^{m-1}(1+\alpha\sigma)^j\big)^{-1} \cdot \big(\sum_{j=0}^{m-1}(1+\alpha\sigma)^j \cdot y_{sm+j+1}\big)$;
**end for**
**return** $\widetilde{x}^S$.

**Algorithm 5** `Katyusha2`$^{\mathsf{ns}}(x_0, S, \sigma, (L_1, \ldots, L_n))$

$m \leftarrow n$; $\overline{L} = (L_1 + \cdots + L_n)/n$;
$\tau_2 \leftarrow \frac{1}{2}$;
$y_0 = z_0 = \widetilde{x}^0 \leftarrow x_0$;
**for** $s \leftarrow 0$ **to** $S-1$ **do**
    $\tau_{1,s} \leftarrow \frac{2}{s+4}$, $\alpha_s \leftarrow \frac{1}{9\tau_{1,s}\overline{L}}$
    $\mu^s \leftarrow \nabla f(\widetilde{x}^s)$;
    **for** $j \leftarrow 0$ **to** $m-1$ **do**
        $k \leftarrow (sm) + j$;
        $x_{k+1} \leftarrow \tau_{1,s} z_k + \tau_2 \widetilde{x}^s + (1 - \tau_{1,s} - \tau_2) y_k$;
        Pick $i$ randomly from $\{1, 2, \ldots, n\}$, each with probability $L_i/n\overline{L}$;
        $\widetilde{\nabla}_{k+1} \leftarrow \mu^s + \nabla f_i(x_{k+1}) - \nabla f_i(\widetilde{x}^s)$;
        $z_{k+1} = \arg\min_z \big\{\frac{1}{\alpha_s} V_{z_k}(z) + \langle \widetilde{\nabla}_{k+1}, z\rangle + \psi(z)\big\}$;
            ◇ *$V_x(y)$ is the Bregman divergence function, see Section 6*
        $y_{k+1} \leftarrow \arg\min_y \big\{\frac{9\overline{L}}{2}\|y - x_{k+1}\|^2 + \langle \widetilde{\nabla}_{k+1}, y\rangle + \psi(y)\big\}$;
    **end for**
    $\widetilde{x}^{s+1} \leftarrow \frac{1}{m}\sum_{j=1}^{m} y_{sm+j}$;
**end for**
**return** $\widetilde{x}^S$.

## E.2 One-Iteration Analysis

Similar as Section 2.1, we first analyze the behavior of `Katyusha2` in a single iteration (i.e., for a fixed $k$). We view $y_k, z_k$ and $x_{k+1}$ as fixed in this section so the only randomness comes from the choice of $i$ in iteration $k$. We abbreviate in this subsection by $\widetilde{x} = \widetilde{x}^s$ where $s$ is the epoch that iteration $k$ belongs to, and denote by $\sigma_{k+1}^2 \stackrel{\text{def}}{=} \|\nabla f(x_{k+1}) - \widetilde{\nabla}_{k+1}\|_*^2$.

Our first lemma is analogous to Lemma E.1 except the change of the parameter and the norm.

**Lemma E.1** (proximal gradient descent). *If*

$$y_{k+1} = \arg\min_y \big\{ \frac{9\overline{L}}{2}\|y - x_{k+1}\|^2 + \langle \widetilde{\nabla}_{k+1}, y - x_{k+1}\rangle + \psi(y) - \psi(x_{k+1}) \big\} \ , \qquad \textit{and}$$

$$\mathsf{Prog}(x_{k+1}) \stackrel{\text{def}}{=} -\min_y \big\{ \frac{9\overline{L}}{2}\|y - x_{k+1}\|^2 + \langle \widetilde{\nabla}_{k+1}, y - x_{k+1}\rangle + \psi(y) - \psi(x_{k+1}) \big\} \geq 0 \ ,$$

*we have (where the expectation is only over the randomness of* $\widetilde{\nabla}_{k+1}$*)*

$$F(x_{k+1}) - \mathbb{E}\big[F(y_{k+1})\big] \geq \mathbb{E}\big[\mathsf{Prog}(x_{k+1})\big] - \frac{1}{16\overline{L}}\mathbb{E}\big[\sigma_{k+1}^2\big] \ .$$

*Proof.*

$$\begin{aligned}
\mathsf{Prog}(x_{k+1}) &= -\min_y \Big\{ \frac{9\overline{L}}{2}\|y - x_{k+1}\|^2 + \langle \widetilde{\nabla}_{k+1}, y - x_{k+1}\rangle + \psi(y) - \psi(x_{k+1}) \Big\} \\
&\overset{①}{=} -\Big( \frac{9\overline{L}}{2}\|y_{k+1} - x_{k+1}\|^2 + \langle \widetilde{\nabla}_{k+1}, y_{k+1} - x_{k+1}\rangle + \psi(y_{k+1}) - \psi(x_{k+1}) \Big) \\
&= -\Big( \frac{\overline{L}}{2}\|y_{k+1} - x_{k+1}\|^2 + \langle \nabla f(x_{k+1}), y_{k+1} - x_{k+1}\rangle + \psi(y_{k+1}) - \psi(x_{k+1}) \Big) \\
&\quad + \Big( \langle \nabla f(x_{k+1}) - \widetilde{\nabla}_{k+1}, y_{k+1} - x_{k+1}\rangle - 4\overline{L}\|y_{k+1} - x_{k+1}\|^2 \Big) \\
&\overset{②}{\leq} -\Big( f(y_{k+1}) - f(x_{k+1}) + \psi(y_{k+1}) - \psi(x_{k+1}) \Big) + \frac{1}{16\overline{L}}\|\nabla f(x_{k+1}) - \widetilde{\nabla}_{k+1}\|_*^2 \ .
\end{aligned}$$

Above, ① is by the definition of $y_{k+1}$, and ② uses the smoothness of function $f(\cdot)$, as well as Young's inequality $\langle a, b\rangle - \frac{1}{2}\|b\|^2 \leq \frac{1}{2}\|a\|_*^2$. Taking expectation on both sides we arrive at the desired result. □

The next lemma is analogous to Lemma 2.4 but with slightly different proof.

**Lemma E.2** (variance upper bound).

$$\mathbb{E}\big[\|\widetilde{\nabla}_{k+1} - \nabla f(x_{k+1})\|_*^2\big] \leq 8\overline{L} \cdot \big(f(\widetilde{x}) - f(x_{k+1}) - \langle \nabla f(x_{k+1}), \widetilde{x} - x_{k+1}\rangle\big) \ .$$

*Proof.* Each $f_i(x)$, being convex and $L_i$-smooth, implies the following inequality which is classical in convex optimization and can be found for instance in Theorem 2.1.5 of the textbook of Nesterov [43].

$$\|\nabla f_i(x_{k+1}) - \nabla f_i(\widetilde{x})\|_*^2 \leq 2L_i \cdot \big(f_i(\widetilde{x}) - f_i(x_{k+1}) - \langle \nabla f_i(x_{k+1}), \widetilde{x} - x_{k+1}\rangle\big) \qquad \text{(E.1)}$$

Therefore, taking expectation over the random choice of $i$, we have

$$\begin{aligned}
&\mathbb{E}\big[\|\widetilde{\nabla}_{k+1} - \nabla f(x_{k+1})\|_*^2\big] \\
&= \mathbb{E}\big[\big\|\frac{1}{np_i}\big(\nabla f_i(x_{k+1}) - \nabla f_i(\widetilde{x})\big) - \big(\nabla f(x_{k+1}) - \nabla f(\widetilde{x})\big)\big\|_*^2\big] \\
&\overset{①}{\leq} 2\mathbb{E}\big[\frac{1}{n^2p_i^2}\big\|\nabla f_i(x_{k+1}) - \nabla f_i(\widetilde{x})\big\|_*^2\big] + 2\big\|\nabla f(x_{k+1}) - \nabla f(\widetilde{x})\big\|_*^2 \\
&\overset{②}{\leq} 4 \cdot \mathbb{E}\big[\frac{L_i}{n^2p_i^2}\big(f_i(\widetilde{x}) - f_i(x_{k+1}) - \langle \nabla f_i(x_{k+1}), \widetilde{x} - x_{k+1}\rangle\big)\big] + 2\big\|\nabla f(x_{k+1}) - \nabla f(\widetilde{x})\big\|_*^2 \\
&\overset{③}{=} 4\overline{L} \cdot \big(f(\widetilde{x}) - f(x_{k+1}) - \langle \nabla f(x_{k+1}), \widetilde{x} - x_{k+1}\rangle\big) + 2\big\|\nabla f(x_{k+1}) - \nabla f(\widetilde{x})\big\|_*^2 \\
&\leq 8\overline{L} \cdot \big(f(\widetilde{x}) - f(x_{k+1}) - \langle \nabla f(x_{k+1}), \widetilde{x} - x_{k+1}\rangle\big) \ .
\end{aligned}$$

Above, inequality ① is because $\|a+b\|_*^2 \le (\|a\|_* + \|b\|_*)^2 \le 2\|a\|_*^2 + 2\|b\|_*^2$; inequality ② follows from (E.1); equality ③ follows from the probability distribution that we select $i$ with probability $p_i = L_i/(n\overline{L})$; inequality ④ uses (E.1) again but replacing $f_i(\cdot)$ with $f(\cdot)$. $\square$

The next lemma is classical for mirror descent with respect to a general Bregman divergence.

**Lemma E.3** (proximal mirror descent)**.** *Suppose $\psi(\cdot)$ is $\sigma$-SC with respect to $V_x(y)$. Then, fixing $\widetilde{\nabla}_{k+1}$ and letting*

$$z_{k+1} = \arg\min_z \left\{ V_{z_k}(z) + \alpha\langle \widetilde{\nabla}_{k+1}, z - z_k\rangle + \alpha\psi(z) - \alpha\psi(z_k) \right\} ,$$

*it satisfies for all $u \in \mathbb{R}^d$,*

$$\alpha\langle \widetilde{\nabla}_{k+1}, z_{k+1} - u\rangle + \alpha\psi(z_{k+1}) - \alpha\psi(u) \le -\frac{1}{2}\|z_k - z_{k+1}\|^2 + V_{z_k}(u) - (1+\alpha\sigma)V_{z_{k+1}}(u) .$$

*Proof.* By the minimality definition of $z_{k+1}$, we have that

$$\nabla V_{z_k}(z_{k+1}) + \alpha\widetilde{\nabla}_{k+1} + \alpha g = 0$$

where $g$ is *some* subgradient of $\psi(z)$ at point $z = z_{k+1}$. This implies that for every $u$ it satisfies

$$0 = \left\langle \nabla V_{z_k}(z_{k+1}) + \alpha\widetilde{\nabla}_{k+1} + \alpha g, z_{k+1} - u \right\rangle .$$

At this point, using the equality $\langle \nabla V_{z_k}(z_{k+1}), z_{k+1} - u\rangle = V_{z_k}(z_{k+1}) - V_{z_k}(u) + V_{z_{k+1}}(u)$ (known as the "three-point equality of Bregman divergence", see [47]), as well as the inequality $\langle g, z_{k+1} - u\rangle \ge \psi(z_{k+1}) - \psi(u) + \sigma V_{z_{k+1}}(u)$ which comes from the strong convexity of $\psi(\cdot)$, we can write

$$\begin{aligned}
&\alpha\langle \widetilde{\nabla}_{k+1}, z_{k+1} - u\rangle + \alpha\psi(z_{k+1}) - \alpha\psi(u) \\
&= -\langle z_{k+1} - z_k, z_{k+1} - u\rangle - \langle \alpha g, z_{k+1} - u\rangle + \alpha\psi(z_{k+1}) - \alpha\psi(u) \\
&\le -V_{z_k}(z_{k+1}) + V_{z_k}(u) - (1+\alpha\sigma)V_{z_{k+1}}(u) .
\end{aligned}$$

Finally, using $V_{z_k}(z_{k+1}) \ge \frac{1}{2}\|z_k - z_{k+1}\|^2$ which comes from the strong convexity of $w(x)$ with respect to $\|\cdot\|$, we complete the proof. $\square$

The following lemma combines Lemma E.1, Lemma E.2 and Lemma E.3 all together, using the special choice of $x_{k+1}$ which is a convex combination of $y_k, z_k$ and $\widetilde{x}$:

**Lemma E.4** (coupling step 1)**.** *If $x_{k+1} = \tau_1 z_k + \tau_2\widetilde{x} + (1-\tau_1-\tau_2)y_k$, where $\tau_1 \le \frac{1}{9\alpha\overline{L}}$ and $\tau_2 = \frac{1}{2}$,*

$$\begin{aligned}
&\alpha\langle \nabla f(x_{k+1}), z_k - u\rangle - \alpha\psi(u) \\
&\le \frac{\alpha}{\tau_1}\Big( F(x_{k+1}) - \mathbb{E}\big[F(y_{k+1})\big] + \tau_2 F(\widetilde{x}) - \tau_2\mathbb{E}\big[F(x_{k+1})\big] - \tau_2\langle \nabla f(x_{k+1}), \widetilde{x} - x_{k+1}\rangle \Big) \\
&\quad + V_{z_k}(u) - (1+\alpha\sigma)\mathbb{E}\big[V_{z_{k+1}}(u)\big] + \frac{\alpha(1-\tau_1-\tau_2)}{\tau_1}\psi(y_k) - \frac{\alpha}{\tau_1}\psi(x_{k+1}) .
\end{aligned}$$

*Proof.* We first apply Lemma E.3 and get

$$\begin{aligned}
&\alpha\langle \widetilde{\nabla}_{k+1}, z_k - u\rangle + \alpha\psi(z_{k+1}) - \alpha\psi(u) \\
&= \alpha\langle \widetilde{\nabla}_{k+1}, z_k - z_{k+1}\rangle + \alpha\langle \widetilde{\nabla}_{k+1}, z_{k+1} - u\rangle + \alpha\psi(z_{k+1}) - \alpha\psi(u) \\
&\le \alpha\langle \widetilde{\nabla}_{k+1}, z_k - z_{k+1}\rangle - \frac{1}{2}\|z_k - z_{k+1}\|^2 + V_{z_k}(u) - (1+\alpha\sigma)V_{z_{k+1}}(u) .
\end{aligned} \tag{E.2}$$

By defining $v \stackrel{\text{def}}{=} \tau_1 z_{k+1} + \tau_2 \widetilde{x} + (1-\tau_1-\tau_2) y_k$, we have $x_{k+1} - v = \tau_1(z_k - z_{k+1})$ and therefore

$$\begin{aligned}
&\mathbb{E}\Big[\alpha\langle\widetilde{\nabla}_{k+1}, z_k - z_{k+1}\rangle - \frac{1}{2}\|z_k - z_{k+1}\|^2\Big] = \mathbb{E}\Big[\frac{\alpha}{\tau_1}\langle\widetilde{\nabla}_{k+1}, x_{k+1} - v\rangle - \frac{1}{2\tau_1^2}\|x_{k+1} - v\|^2\Big] \\
&= \mathbb{E}\Big[\frac{\alpha}{\tau_1}\Big(\langle\widetilde{\nabla}_{k+1}, x_{k+1} - v\rangle - \frac{1}{2\alpha\tau_1}\|x_{k+1} - v\|^2 - \psi(v) + \psi(x_{k+1})\Big) + \frac{\alpha}{\tau_1}\Big(\psi(v) - \psi(x_{k+1})\Big)\Big] \\
&\overset{①}{\leq} \mathbb{E}\Big[\frac{\alpha}{\tau_1}\Big(\langle\widetilde{\nabla}_{k+1}, x_{k+1} - v\rangle - \frac{9\overline{L}}{2}\|x_{k+1} - v\|^2 - \psi(v) + \psi(x_{k+1})\Big) + \frac{\alpha}{\tau_1}\Big(\psi(v) - \psi(x_{k+1})\Big)\Big] \\
&\overset{②}{\leq} \mathbb{E}\Big[\frac{\alpha}{\tau_1}\Big(F(x_{k+1}) - F(y_{k+1}) + \frac{1}{16\overline{L}}\sigma_{k+1}^2\Big) + \frac{\alpha}{\tau_1}\Big(\psi(v) - \psi(x_{k+1})\Big)\Big] \\
&\overset{③}{\leq} \mathbb{E}\Big[\frac{\alpha}{\tau_1}\Big(F(x_{k+1}) - F(y_{k+1}) + \frac{1}{2}\big(f(\widetilde{x}) - f(x_{k+1}) - \langle\nabla f(x_{k+1}), \widetilde{x} - x_{k+1}\rangle\big)\Big) \\
&\quad + \frac{\alpha}{\tau_1}\Big(\tau_1\psi(z_{k+1}) + \tau_2\psi(\widetilde{x}) + (1-\tau_1-\tau_2)\psi(y_k) - \psi(x_{k+1})\Big)\Big] \ . \qquad \text{(E.3)}
\end{aligned}$$

Above, ① uses our choice $\tau_1 \leq \frac{1}{9\alpha\overline{L}}$, ② uses Lemma E.1, ③ uses Lemma E.2 together with the convexity of $\psi(\cdot)$ and the definition of $v$. Finally, noticing that $\mathbb{E}[\langle\widetilde{\nabla}_{k+1}, z_k - u\rangle] = \langle\nabla f(x_{k+1}), z_k - u\rangle$ and $\tau_2 = \frac{1}{2}$, we obtain the desired inequality by combining (E.2) and (E.3). □

The next lemma is completely analogous to Lemma 2.7 except that we use Lemma E.4 rather than Lemma 2.6. We ignore the proof since it is a simple copy-and-paste.

**Lemma E.5** (coupling step 2)**.** *Under the same choices of $\tau_1, \tau_2$ as in Lemma E.4, we have*

$$\begin{aligned}
0 \leq \frac{\alpha(1-\tau_1-\tau_2)}{\tau_1}(F(y_k) - F(x^*)) - \frac{\alpha}{\tau_1}\big(\mathbb{E}\big[F(y_{k+1})\big] - F(x^*)\big) + \frac{\alpha\tau_2}{\tau_1}\big(F(\widetilde{x}) - \tau_2 F(x^*)\big) \\
+ V_{z_k}(x^*) - (1+\alpha\sigma)\mathbb{E}\big[V_{z_{k+1}}(x^*)\big] \ .
\end{aligned}$$

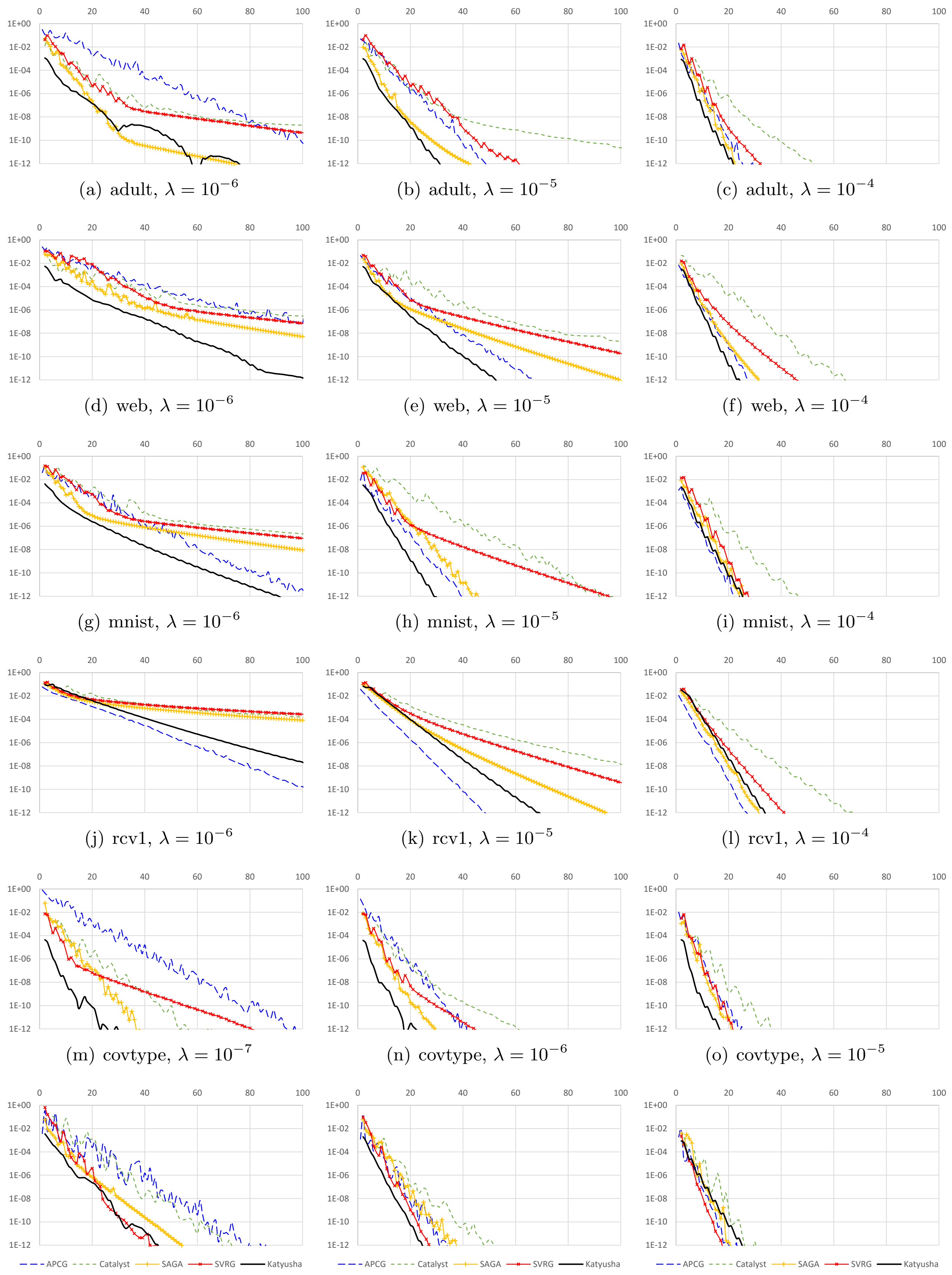

(a) adult, $\lambda = 10^{-6}$ (b) adult, $\lambda = 10^{-5}$ (c) adult, $\lambda = 10^{-4}$

(d) web, $\lambda = 10^{-6}$ (e) web, $\lambda = 10^{-5}$ (f) web, $\lambda = 10^{-4}$

(g) mnist, $\lambda = 10^{-6}$ (h) mnist, $\lambda = 10^{-5}$ (i) mnist, $\lambda = 10^{-4}$

(j) rcv1, $\lambda = 10^{-6}$ (k) rcv1, $\lambda = 10^{-5}$ (l) rcv1, $\lambda = 10^{-4}$

(m) covtype, $\lambda = 10^{-7}$ (n) covtype, $\lambda = 10^{-6}$ (o) covtype, $\lambda = 10^{-5}$

(p) sensit, $\lambda = 10^{-6}$ (q) sensit, $\lambda = 10^{-5}$ (r) sensit, $\lambda = 10^{-4}$

Figure 3: Experiments on ridge regression with $\ell_2$ regularizer weight $\lambda$.

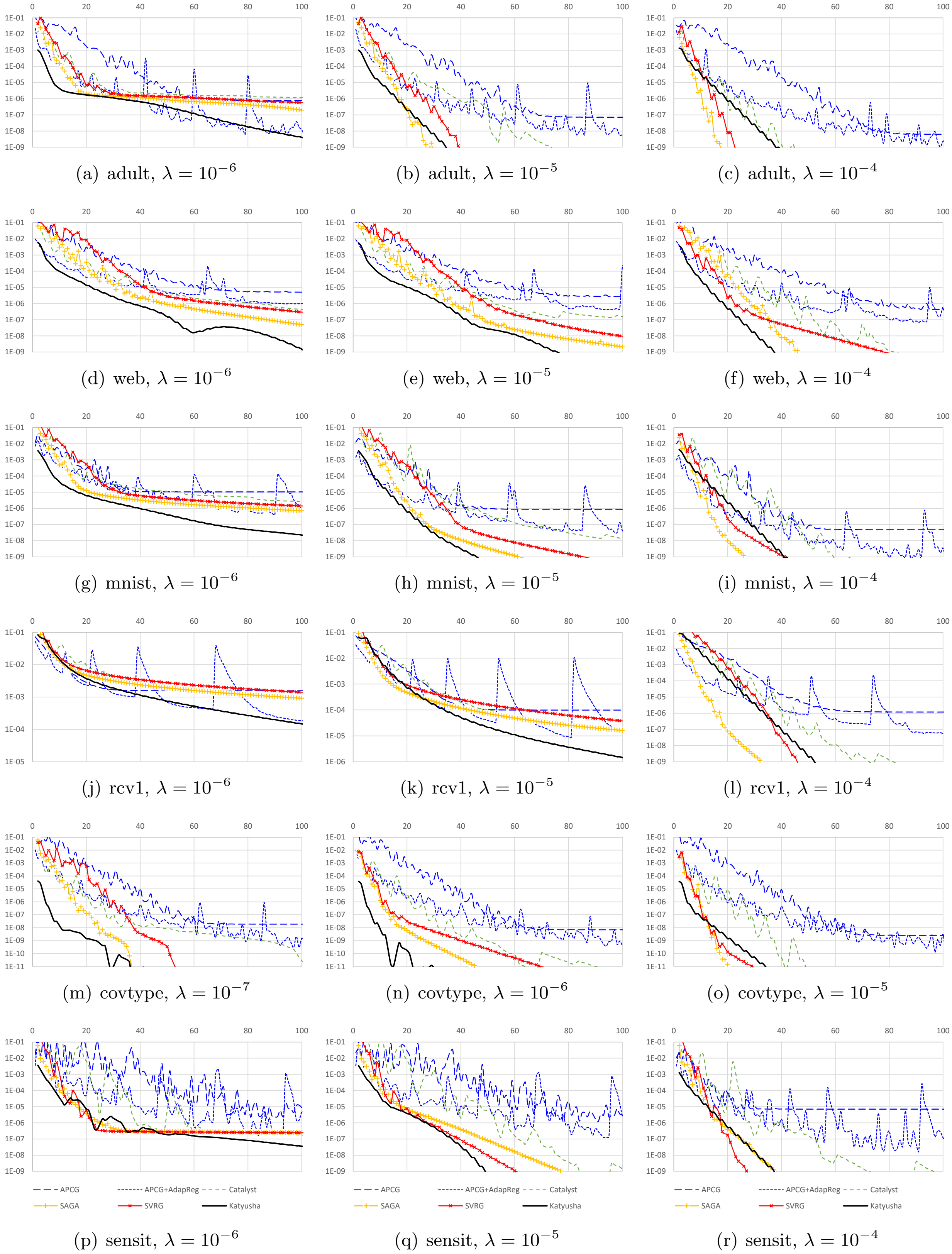

(a) adult, $\lambda = 10^{-6}$ (b) adult, $\lambda = 10^{-5}$ (c) adult, $\lambda = 10^{-4}$

(d) web, $\lambda = 10^{-6}$ (e) web, $\lambda = 10^{-5}$ (f) web, $\lambda = 10^{-4}$

(g) mnist, $\lambda = 10^{-6}$ (h) mnist, $\lambda = 10^{-5}$ (i) mnist, $\lambda = 10^{-4}$

(j) rcv1, $\lambda = 10^{-6}$ (k) rcv1, $\lambda = 10^{-5}$ (l) rcv1, $\lambda = 10^{-4}$

(m) covtype, $\lambda = 10^{-7}$ (n) covtype, $\lambda = 10^{-6}$ (o) covtype, $\lambda = 10^{-5}$

(p) sensit, $\lambda = 10^{-6}$ (q) sensit, $\lambda = 10^{-5}$ (r) sensit, $\lambda = 10^{-4}$

Figure 4: Experiments on Lasso with $\ell_1$ regularizer weight $\lambda$.